# Impacto del Enfoque "Matemáticas en Tres Actos" en la Educación Matemática

# Impact of the "Mathematics in Three Acts" Approach on Mathematics Education

**Resumen:**

El enfoque "Matemáticas en Tres Actos", propuesto por Dan Meyer, busca transformar la enseñanza de las matemáticas a través de un modelo que fomenta la participación activa de los estudiantes, promoviendo la creatividad, la resolución de problemas y la metacognición. Este estudio explora la implementación de este enfoque en un contexto de concurso matemático para estudiantes de secundaria, evaluando su impacto en diversas competencias clave. Se examinan aspectos como la creatividad matemática, la capacidad para resolver problemas, las habilidades metacognitivas y las percepciones de los estudiantes hacia las matemáticas. Los resultados muestran que el enfoque contribuye al desarrollo de habilidades creativas, mejora la comprensión y resolución de problemas, y aumenta la motivación y la confianza de los estudiantes. Sin embargo, también se identifican áreas de mejora, especialmente en la justificación de procedimientos y en la flexibilidad cognitiva. Este estudio resalta la efectividad del enfoque "Matemáticas en Tres Actos" como una metodología innovadora que favorece un aprendizaje más significativo, reflexivo y autónomo, sugiriendo su potencial para transformar la enseñanza de las matemáticas en contextos educativos diversos.

**Palabras clave:** Matemáticas en Tres Actos, educación matemática, creatividad matemática, resolución de problemas, habilidades metacognitivas, motivación estudiantil, enseñanza activa.

**Abstract:**

The "Mathematics in Three Acts" approach, proposed by Dan Meyer, aims to transform the teaching of mathematics through a model that encourages active student participation, fostering creativity, problem-solving, and metacognition. This study explores the implementation of this approach in a mathematics contest for secondary school students, evaluating its impact on various key competencies. Aspects such as mathematical creativity, problem-solving skills, metacognitive abilities, and students' perceptions of mathematics are examined. The results show that the approach contributes to the development of creative skills, improves understanding and problem-solving abilities, and increases student motivation and confidence. However, areas for improvement are also identified, particularly in the justification of procedures and cognitive flexibility. This study highlights the effectiveness of the "Mathematics in Three Acts" approach as an innovative methodology that fosters more meaningful, reflective, and autonomous learning, suggesting its potential to transform mathematics teaching in diverse educational contexts.

**Resumo:**

O enfoque "Matemática em Três Atos", proposto por Dan Meyer, busca transformar o ensino de matemática por meio de um modelo que incentiva a participação ativa dos estudantes, promovendo a criatividade, a resolução de problemas e a metacognição. Este estudo explora a implementação desse enfoque em um concurso de matemática para alunos do ensino secundário, avaliando seu impacto em várias competências-chave. São analisados aspectos como a criatividade matemática, a capacidade de resolver problemas, as habilidades metacognitivas e as percepções dos estudantes sobre a matemática. Os resultados mostram que o enfoque contribui para o desenvolvimento das habilidades criativas, melhora a compreensão e a resolução de problemas, e aumenta a motivação e a confiança dos estudantes. No entanto, também foram identificadas áreas de melhoria, especialmente na justificativa dos procedimentos e na flexibilidade cognitiva. Este estudo destaca a eficácia do enfoque "Matemática em Três Atos" como uma metodologia inovadora que favorece um aprendizado mais significativo, reflexivo e autônomo, sugerindo seu potencial para transformar o ensino de matemática em contextos educacionais diversos.

**Palavras chave:** Matemática em Três Atos, educação matemática, criatividade matemática, resolução de problemas, habilidades metacognitivas, motivação estudantil, ensino ativo.

**Résumé:**

L'approche "Mathématiques en Trois Actes", proposée par Dan Meyer, vise à transformer l'enseignement des mathématiques par un modèle encourageant la participation active des étudiants, favorisant la créativité, la résolution de problèmes et la métacognition. Cette étude explore la mise en œuvre de cette approche dans un concours de mathématiques pour les élèves du secondaire, en évaluant son impact sur diverses compétences clés. Les aspects tels que la créativité mathématique, la résolution de problèmes, les compétences métacognitives et les perceptions des élèves vis-à-vis des mathématiques sont examinés. Les résultats montrent que l'approche contribue au développement des compétences créatives, améliore la compréhension et la résolution de problèmes, et accroît la motivation et la confiance des étudiants. Toutefois, des domaines d'amélioration sont également identifiés, notamment dans la justification des procédures et la flexibilité cognitive. Cette étude souligne l'efficacité de l'approche "Mathématiques en Trois Actes" en tant que méthodologie innovante qui favorise un apprentissage plus significatif, réfléchi et autonome, suggérant son potentiel pour transformer l'enseignement des mathématiques dans divers contextes éducatifs.

**Mots clés:** Mathématiques en Trois Actes, éducation mathématique, créativité mathématique, résolution de problèmes, compétences métacognitives, motivation des élèves, enseignement actif.

# 1. Introducción

Las matemáticas no solo proporcionan conocimientos técnicos, sino que también fomentan habilidades esenciales como el pensamiento lógico, crítico y creativo, las cuales son fundamentales para enfrentar problemas complejos y tomar decisiones informadas en diversas áreas de la vida (Ayllón, Gómez, y Ballesta-Claver, 2016). Sin embargo, la enseñanza de esta disciplina enfrenta múltiples desafíos, particularmente cuando se aplican métodos tradicionales centrados en la memorización y la repetición de procedimientos. Estos enfoques limitan una comprensión profunda de los conceptos matemáticos y disminuyen la motivación de los estudiantes, quienes a menudo perciben las matemáticas como abstractas y carentes de aplicación práctica (Valero y González, 2020).

Para superar estas limitaciones, es imprescindible adoptar enfoques pedagógicos innovadores que promuevan un aprendizaje más participativo y significativo. Un enfoque destacado es el de Matemáticas en Tres Actos, propuesto por Dan Meyer. Este modelo narrativo, que se estructura en tres fases, permite a los estudiantes analizar problemas y explorar de manera autónoma los datos necesarios para resolverlos, lo cual estimula su creatividad, autonomía y capacidad crítica (Meyer, 2011). Al involucrar activamente a los estudiantes en la formulación y resolución de problemas, este enfoque favorece una conexión más profunda con los conceptos matemáticos y desarrolla competencias clave como la formulación de hipótesis, la reflexión y la evaluación crítica de soluciones (Cruz, Arteaga y Del Sol Martínez, 2018).

La creatividad matemática, junto con la capacidad para resolver problemas y el desarrollo de habilidades metacognitivas, resulta ser especialmente relevante en este enfoque. Al permitir que los estudiantes enfrenten problemas desde diferentes perspectivas, formulen preguntas y busquen soluciones no convencionales, se promueve un aprendizaje más profundo. Estudios previos han demostrado que la creatividad en las matemáticas contribuye al fortalecimiento de la confianza de los estudiantes y reduce su ansiedad hacia la asignatura (Ibeas, 2020; Leikin y Pitta-Pantazi, 2013). Además, un enfoque que favorezca el descubrimiento activo y la experimentación no solo mejora el rendimiento académico, sino que también incrementa la motivación hacia las matemáticas, creando un ambiente en el que los estudiantes se sienten comprometidos con su aprendizaje (Boaler, 2016; Valero y González, 2020).

Este estudio es relevante ya que evalúa la implementación de una adaptación del enfoque Matemáticas en Tres Actos en el contexto de un concurso de matemáticas para estudiantes de secundaria. Más allá de intentar superar las limitaciones de los métodos tradicionales, este estudio busca potenciar habilidades como la creatividad, la motivación y la comprensión matemática, con el objetivo de fomentar un aprendizaje más significativo y efectivo.

La enseñanza de las matemáticas en secundaria enfrenta desafíos significativos en cuanto a la promoción de la creatividad, el desarrollo de habilidades metacognitivas y la resolución de problemas. A pesar de la importancia de estas competencias para el pensamiento crítico y la formación integral, los métodos tradicionales de enseñanza siguen centrados en la memorización y la aplicación mecánica de procedimientos. Esta rigidez en la enseñanza no solo limita el aprendizaje profundo, sino que también afecta la motivación de los estudiantes y su percepción de las matemáticas como una disciplina relevante y aplicable a su vida cotidiana (Valero y González, 2020). Como consecuencia, muchos estudiantes ven las matemáticas como una materia abstracta, sin conexión con la realidad, lo que reduce su disposición a participar activamente y explorar soluciones novedosas (Cruz, Arteaga y Del Sol Martínez, 2018).

Además, hay una falta de metodologías que fomenten la autonomía, el pensamiento crítico y una actitud positiva hacia las matemáticas. El enfoque de Matemáticas en Tres Actos representa una alternativa innovadora que ha mostrado resultados positivos en diversos contextos educativos internacionales, al promover la participación activa y la exploración creativa de problemas (Meyer, 2011). Sin embargo, son limitados los estudios que investigan su impacto en contextos hispanohablantes, especialmente en actividades de matemáticas aplicadas a concursos a nivel secundario. Esta carencia de investigaciones impide comprender plenamente cómo este enfoque puede contribuir al desarrollo de habilidades creativas y metacognitivas, así como mejorar las percepciones de los estudiantes hacia las matemáticas en el ámbito educativo de habla hispana.

Por tanto, el problema de investigación que este estudio pretende abordar es evaluar el impacto de la adaptación del enfoque Matemáticas en Tres Actos en la creatividad matemática, la resolución de problemas, las habilidades metacognitivas y las percepciones de los estudiantes hacia las matemáticas. Este estudio busca analizar si este enfoque innovador puede ofrecer una alternativa efectiva a los métodos tradicionales, promoviendo un aprendizaje significativo que permita a los estudiantes formular, resolver y reflexionar sobre problemas matemáticos, mientras se fomenta una percepción positiva de la disciplina en un contexto desafiante como un concurso académico.

**Objetivo General**

Evaluar el impacto de la adaptación del enfoque de Matemáticas en Tres Actos en la creatividad matemática, la resolución de problemas, las habilidades metacognitivas y las percepciones hacia las matemáticas en estudiantes de secundaria.

**Objetivos Específicos**

1. Evaluar la capacidad de los estudiantes para generar y formular problemas matemáticos de manera creativa.
2. Analizar las competencias de los estudiantes en la resolución de problemas matemáticos.
3. Explorar el desarrollo de habilidades metacognitivas durante el proceso de creación y resolución de problemas.
4. Investigar las percepciones de los estudiantes sobre el impacto del enfoque en su motivación, actitud y confianza hacia las matemáticas.

Desde una perspectiva teórica, esta investigación busca enriquecer el conocimiento existente sobre metodologías innovadoras en la enseñanza de las matemáticas. Diversos estudios han resaltado la necesidad de adoptar enfoques educativos que vayan más allá de la memorización y la aplicación mecánica de procedimientos, favoreciendo la comprensión profunda, la creatividad y el pensamiento crítico (Ayllón et al., 2016; Boaler, 2016). Sin embargo, la investigación educativa aún carece de estudios que profundicen en cómo estas metodologías pueden ser adaptadas y aplicadas en contextos específicos, como los de habla hispana, donde las prácticas tradicionales siguen predominando (Valero y González, 2020).

El enfoque de Matemáticas en Tres Actos, desarrollado por Dan Meyer, se destaca por involucrar a los estudiantes en un proceso narrativo que les permite formular y resolver problemas matemáticos de manera activa (Meyer, 2011). Este estudio tiene como propósito adaptar esta metodología al contexto de un concurso matemático en educación secundaria y evaluar su impacto en el desarrollo de habilidades clave como la creatividad, la metacognición y la autonomía en la resolución de problemas.

Desde una perspectiva práctica, este estudio tiene el potencial de mejorar las prácticas docentes y transformar las estrategias de enseñanza en el aula de matemáticas. Proporcionando un análisis detallado de la implementación del enfoque de Matemáticas en Tres Actos, el estudio ofrece a los docentes un modelo replicable y adaptable que puede integrarse en diversos contextos educativos, promoviendo un aprendizaje activo donde los estudiantes asuman un rol protagónico en su proceso de aprendizaje (Valero y González, 2020).

Este enfoque adaptado no solo fomenta la motivación y el compromiso de los estudiantes con las matemáticas, sino que también les permite enfrentar problemas abiertos que favorecen el desarrollo de habilidades matemáticas y metacognitivas, a la vez que fortalece su confianza para abordar desafíos complejos (Ibeas, 2020). De este modo, el estudio responde a la creciente necesidad de metodologías que preparen a los estudiantes para el mundo real, donde la creatividad, el pensamiento crítico y la capacidad para resolver problemas son esenciales (OCDE, 2019).

# 2. Marco Teórico

El enfoque adaptado "Matemáticas en Tres Actos" se fundamenta en principios pedagógicos clave que buscan transformar la enseñanza de las matemáticas, haciendo el aprendizaje más significativo, activo y centrado en el estudiante. Este marco teórico explora las teorías y conceptos que sustentan dicha adaptación, incluyendo el constructivismo, el aprendizaje basado en problemas (ABP), la creatividad matemática, la resolución de problemas y la metacognición.

## 2.1 El Constructivismo y la Educación Matemática

El constructivismo, propuesto por Piaget (1970), es uno de los pilares fundamentales sobre los que se basa el enfoque "Matemáticas en Tres Actos". Según Piaget, el conocimiento no se transmite de manera pasiva, sino que se construye activamente por el estudiante mientras interactúa con su entorno y resuelve problemas. Este proceso implica tanto la asimilación (incorporación de nueva información en estructuras cognitivas previas) como la acomodación (modificación de esas estructuras para adaptarse a nueva información).

En el enfoque adaptado "Matemáticas en Tres Actos", se promueve un aprendizaje activo mediante la creación y resolución de problemas contextualizados, lo que permite a los estudiantes construir su conocimiento a través de la exploración, formulación de problemas y reflexión sobre sus soluciones. Este enfoque facilita un aprendizaje profundo y significativo que va más allá de la simple memorización de procedimientos, y fomenta la capacidad de los estudiantes para hacer conexiones significativas entre los conceptos matemáticos.

## 2.2 La Teoría Sociocultural de Vygotsky

Complementando el constructivismo, Vygotsky (1978) ofrece una perspectiva sociocultural que destaca la importancia de la interacción social en el proceso de aprendizaje. Según Vygotsky, el conocimiento se desarrolla dentro de un contexto social y cultural, principalmente a través de la interacción con otros. Un concepto clave en esta teoría es la "zona de desarrollo próximo" (ZDP), que se refiere a la distancia entre lo que un estudiante puede hacer de forma independiente y lo que puede hacer con la ayuda de un compañero más capacitado o un docente.

En el enfoque adaptado "Matemáticas en Tres Actos", se promueve el aprendizaje colaborativo al permitir que los estudiantes trabajen en grupos para resolver problemas y generar nuevas ideas. La colaboración y la justificación de las soluciones no solo contribuyen a la construcción del conocimiento individual, sino que también facilitan la construcción colectiva. Este ambiente colaborativo fomenta un contexto de aprendizaje dinámico, donde los estudiantes se ayudan mutuamente a superar dificultades y desarrollan habilidades de comunicación matemática.

## 2.3 Aprendizaje Basado en Problemas (ABP)

El Aprendizaje Basado en Problemas (ABP) se presenta como un enfoque pedagógico que coloca a los estudiantes en el centro del proceso de aprendizaje. En lugar de centrarse únicamente en la transmisión de información, el ABP fomenta la resolución de problemas auténticos y desafiantes, estimulando habilidades críticas como el pensamiento lógico, la creatividad y la resolución autónoma de problemas. El enfoque "Matemáticas en Tres Actos" comparte principios con el ABP al presentar problemas contextualizados que no tienen soluciones predefinidas, lo que requiere pensamiento divergente y justificación de soluciones.

De acuerdo con Hiebert y Carpenter (1992), el ABP permite a los estudiantes desarrollar una comprensión profunda de los conceptos matemáticos, ya que se enfrentan a problemas que requieren la aplicación de diversos procedimientos y herramientas matemáticas. En el enfoque adaptado "Matemáticas en Tres Actos", la fase de exploración (Acto 1) permite que los estudiantes formulen preguntas y piensen en posibles soluciones, lo que fomenta la investigación activa y el desarrollo de habilidades para resolver problemas en contextos reales.

## 2.4 Creatividad Matemática

La creatividad matemática es un concepto esencial dentro del enfoque adaptado "Matemáticas en Tres Actos", ya que no solo permite a los estudiantes resolver problemas de manera original, sino también formular nuevos problemas y hacer conexiones entre distintas áreas de las matemáticas. Silver (1997) sostiene que la creatividad matemática debe ser entendida como un proceso abierto, donde los estudiantes exploran múltiples soluciones posibles y reflexionan sobre el proceso de resolución.

Esta variante del enfoque fomenta la creatividad especialmente en la fase de creación del problema (Acto 2), donde los estudiantes tienen la libertad de generar sus propios problemas. Esto les otorga mayor autonomía y control sobre su aprendizaje, estimulando el pensamiento flexible y divergente sobre los problemas presentados. Al involucrar a los estudiantes en la creación y formulación de problemas, el enfoque "Matemáticas en Tres Actos" facilita la generación de soluciones creativas y el desarrollo de habilidades de pensamiento crítico.

Para evaluar la creatividad matemática en este estudio, se definieron varios indicadores que permiten medir diferentes aspectos del proceso creativo de los estudiantes. Estos indicadores incluyen la diversidad y fluidez en la formulación de preguntas, la capacidad para reinterpretar contextos de manera original y la innovación en el enfoque del problema. Además, se mide la pertinencia de las preguntas formuladas, la calidad de los datos proporcionados en el problema y la existencia de un entorno matemático coherente.

## 2.5 Resolución de Problemas Matemáticos

La resolución de problemas matemáticos es una habilidad central en la educación matemática, ya que no solo implica aplicar conceptos y procedimientos para encontrar soluciones, sino también desarrollar un pensamiento crítico y analítico (Polya, 1945). La capacidad para resolver problemas permite a los estudiantes enfrentarse a situaciones no rutinarias y adaptar sus estrategias según el contexto (Lester, 1994). Esta habilidad es crucial no solo en el aula, sino también en la vida cotidiana y en diversas situaciones reales donde se requiere flexibilidad y análisis crítico (OECD, 2014).

Para evaluar las competencias de los estudiantes en la resolución de problemas, se utilizan indicadores como la capacidad para comunicar la comprensión del problema, traducir problemas a un lenguaje matemático adecuado y la pertinencia de las estrategias y procedimientos utilizados para la solución. Además, se evalúa la justificación de los procedimientos y los resultados obtenidos, lo que permite medir la capacidad de los estudiantes para argumentar y respaldar sus soluciones con base en razonamientos sólidos.

La implementación del enfoque adaptado "Matemáticas en Tres Actos" se alinea perfectamente con el desarrollo de estas competencias. En el Acto 2, los estudiantes tienen la oportunidad de diseñar y ejecutar planes de solución, lo que les permite enfrentar problemas complejos, aplicar sus conocimientos y ajustar sus estrategias según sea necesario. La reflexión y discusión en el Acto 3 refuerzan este proceso, ayudando a consolidar las competencias de resolución de problemas de los estudiantes.

## 2.6 Metacognición y Autorregulación del Aprendizaje

La metacognición, definida por Flavell (1979), es la capacidad de los estudiantes para reflexionar sobre su propio proceso de pensamiento y regular sus estrategias de aprendizaje. Este proceso incluye tanto el conocimiento metacognitivo (la comprensión de los propios procesos cognitivos) como la regulación metacognitiva (la capacidad de controlar y ajustar estos procesos para mejorar el aprendizaje).

En este estudio, la evaluación de las habilidades metacognitivas se realiza a través de indicadores como la autonomía en la toma de decisiones, la flexibilidad en la estructuración del problema, y la capacidad para identificar y corregir errores. Además, se mide la autorregulación y el ajuste de estrategias durante el proceso de resolución de problemas.

El enfoque adaptado "Matemáticas en Tres Actos" favorece el desarrollo de estas habilidades metacognitivas al permitir que los estudiantes planifiquen su enfoque (Acto 1), monitoreen su comprensión y ajusten sus estrategias (Acto 2), y reflexionen sobre su proceso y resultados (Acto 3). Este proceso promueve la autorregulación y la reflexión crítica, elementos esenciales para el desarrollo metacognitivo de los estudiantes.

## 2.7 Implicaciones para la Educación Matemática

El enfoque adaptado "Matemáticas en Tres Actos" tiene importantes implicaciones para la enseñanza de las matemáticas. Al promover un aprendizaje activo, colaborativo y reflexivo, este enfoque no solo desarrolla habilidades de resolución de problemas, creatividad matemática y autonomía, sino que también integra la metacognición en el proceso educativo. De esta manera, va más allá de la simple memorización de procedimientos, permitiendo a los estudiantes desarrollar una comprensión profunda y significativa de los conceptos matemáticos. La implementación de este enfoque tiene el potencial de transformar la enseñanza de las matemáticas, haciéndola más relevante, dinámica y centrada en el estudiante.

# 3. Metodología

La metodología es un componente esencial de este estudio, ya que describe el diseño de investigación, los participantes, los procedimientos, los instrumentos de recolección de datos y el análisis de los mismos. Una presentación clara y detallada de esta sección asegura la replicabilidad del estudio y permite a los lectores comprender el proceso seguido para obtener las conclusiones.

## 3.1. Diseño del Estudio

Este estudio adopta un enfoque mixto que combina métodos cualitativos y cuantitativos para analizar integralmente el impacto de la adaptación del enfoque "Matemáticas en Tres Actos" en un concurso grupal de matemáticas. El componente cualitativo se apoya en un diseño de estudio de caso, lo que permite explorar en profundidad las experiencias, percepciones y dinámicas grupales de los estudiantes al aplicar el enfoque adaptado. Este diseño es particularmente adecuado para observar cómo emergen la creatividad matemática, la resolución de problemas y las habilidades metacognitivas en un contexto real y desafiante (Stake, 2020).

La recolección de datos cualitativos se realizó mediante observación participante sistemática no controlada (Supo & Zacarías, 2020). La implementación del enfoque "Matemáticas en Tres Actos" requirió ajustarlo de su formato original, concebido para entornos de aula individuales, a un ambiente competitivo que demanda colaboración, creatividad y resolución de problemas dentro de un tiempo limitado. El concurso, dirigido a estudiantes de secundaria organizados en grupos heterogéneos, estructuró sus actividades en torno a las tres fases principales del enfoque: Exploración, Creación y Solución. Este escenario auténtico proporcionó un contexto idóneo para evaluar el impacto del enfoque en las competencias matemáticas de los participantes.

El componente cuantitativo se incorpora mediante el análisis de datos recopilados a través de encuestas aplicadas al finalizar el concurso, utilizando escalas Likert como instrumento. Estas encuestas recogen información sobre las percepciones y actitudes de los estudiantes hacia el enfoque adaptado, midiendo aspectos como la motivación, la confianza y la actitud hacia las matemáticas. La inclusión de este componente cuantitativo complementa los hallazgos cualitativos, permitiendo una triangulación de datos y fortaleciendo la validez de las conclusiones (Creswell & Creswell, 2018).

La sinergia metodológica entre los enfoques cualitativo y cuantitativo es esencial para cumplir con los objetivos de la investigación, ya que permite analizar cómo la adaptación del enfoque "Matemáticas en Tres Actos" influye en diferentes dimensiones del aprendizaje matemático y en las percepciones de los estudiantes. Al centrarse en un evento específico como el concurso de matemáticas, el estudio enriquece la comprensión de cómo las metodologías activas pueden adaptarse y generar resultados significativos en ambientes no convencionales. Así, el enfoque mixto con estudio de caso se revela como una estrategia metodológica idónea para investigar fenómenos complejos en su contexto real (Yin, 2018), proporcionando una visión integral y contextualizada del impacto de la metodología innovadora implementada.

## 3.2. Adaptación del Enfoque Matemáticas en Tres Actos

El presente estudio requirió adaptar el enfoque Matemáticas en Tres Actos para ajustarlo a la naturaleza competitiva y colaborativa de un concurso grupal de matemáticas. Originalmente diseñado para aplicarse de manera individual en entornos de aula con tiempos más flexibles (Meyer, 2011), el enfoque fue modificado para alinearse con las dinámicas grupales y el tiempo limitado del concurso, sin perder su énfasis central en la exploración, creación y resolución de problemas.

La justificación de esta adaptación radica en la necesidad de fomentar la resolución colaborativa de problemas en un contexto grupal, acomodar los tiempos restringidos del concurso—una sesión de dos horas—y cumplir con los objetivos de promover la creatividad, el pensamiento crítico y las habilidades metacognitivas entre los participantes. Además, se buscó mantener la relevancia y el interés contextual de los problemas presentados, aspectos fundamentales en el enfoque original.

La estructura de la adaptación conservó las tres fases principales del enfoque, pero fue modificada para ajustarse a las dinámicas del concurso:

- Acto 1: Exploración. Los estudiantes observaron un escenario visual ya sea imagen o video, diseñado para despertar curiosidad y generar preguntas matemáticas. Cada grupo formuló y discutió diversas preguntas, seleccionando una para desarrollarla en las fases posteriores. Este proceso estimuló la creatividad y la diversidad de perspectivas dentro de los equipos.
- Acto 2: Creación de Problemas. Los grupos trabajaron colectivamente para contextualizar la pregunta seleccionada en un problema matemático coherente. Identificaron los datos necesarios, definieron variables y redactaron un enunciado claro, promoviendo la originalidad y la capacidad de contextualizar conceptos matemáticos de manera precisa.

- Acto 3: Solución. Los grupos resolvieron el problema creado aplicando conceptos y razonamientos matemáticos adecuados. Reflexionaron sobre el proceso de resolución, evaluando la precisión, eficiencia y coherencia de sus soluciones. Este acto desarrolló el pensamiento crítico, permitió justificar los métodos utilizados y fomentó la reflexión metacognitiva entre los participantes.

La implementación de esta adaptación se llevó a cabo en un formato estructurado pero flexible. Cada grupo recibió un cuadernillo de trabajo que incluía el escenario visual para el Acto 1 y las indicaciones para la formulación de preguntas, creación y resolución de problemas en las fases subsiguientes. Los facilitadores, el investigador principal y docentes colaboradores brindaron apoyo logístico y observaron las interacciones grupales, sin intervenir en las decisiones de los estudiantes. El diseño buscó simular un entorno competitivo auténtico, manteniendo el enfoque educativo y propiciando un ambiente propicio para el aprendizaje activo.

### 3.3. Participantes

El estudio incluyó a 55 estudiantes de secundaria que participaron en un concurso de matemáticas organizado por la Institución Educativa. Los participantes, con edades entre 13 y 17 años, provenían de todos los niveles de secundaria, lo que permitió observar las respuestas de distintos grupos al enfoque adaptado en un ambiente competitivo y colaborativo, esencial para los objetivos de la investigación.

En términos de género, el grupo estuvo compuesto por 24 mujeres y 31 hombres, lo que posibilitó un análisis equilibrado de las posibles variaciones en la experiencia y percepción del enfoque según el género. Los estudiantes se organizaron en 18 grupos heterogéneos para facilitar el trabajo en equipo y equilibrar las habilidades. Cada grupo fue identificado con un código para facilitar el análisis de los datos recopilados.

La distribución de los grupos por grado fue la siguiente:
- Primero de secundaria: 4 grupos (P-01, P-02, P-03, P-04)
- Segundo de secundaria: 4 grupos (S-01, S-02, S-03, S-04)
- Tercero de secundaria: 3 grupos (T-01, T-02, T-03)
- Cuarto de secundaria: 5 grupos (C-01, C-02, C-03, C-04, C-05)
- Quinto de secundaria: 2 grupos (Q-01, Q-02)

Esta organización fue realizada por una de las docentes de matemáticas con el fin de equilibrar las habilidades y garantizar una representación equitativa de los distintos grados en cada grupo.

Todos los estudiantes de secundaria de la institución participaron en el estudio, sin aplicar criterios de exclusión, lo que permitió evaluar el impacto del enfoque adaptado en un contexto educativo real y representativo. Este enfoque está alineado con el diseño exploratorio-descriptivo y observacional de la metodología empleada.

El concurso se llevó a cabo en una institución educativa privada urbana, que dispone de recursos y personal docente capacitado para apoyar el uso de metodologías innovadoras en matemáticas. Con un total de 55 estudiantes de secundaria y tres docentes de matemáticas, la institución facilita la implementación de enfoques innovadores, como el concurso de matemáticas. Los docentes participaron activamente como moderadores y colaboradores, brindando apoyo y supervisión a los grupos, así como contribuyendo a la recolección de datos cualitativos mediante la observación participante.

Los estudiantes están familiarizados con un enfoque pedagógico que promueve el aprendizaje activo y colaborativo, lo que favoreció la adopción del enfoque adaptado en el contexto del concurso. Este ambiente educativo facilitó la implementación exitosa del enfoque Matemáticas en Tres Actos, proporcionando un escenario auténtico para evaluar su impacto en las competencias matemáticas y en las percepciones de los estudiantes.

## 3.4. Procedimiento

El procedimiento para implementar la adaptación del enfoque Matemáticas en Tres Actos en el contexto del concurso de matemáticas se diseñó para favorecer tanto el trabajo grupal como la estructura del evento. Esta implementación permitió observar cómo los estudiantes aplicaban sus conocimientos previos, generaban y formulaban problemas matemáticos de forma creativa, resolvían dichos problemas, desarrollaban habilidades metacognitivas y experimentaban cambios en su motivación y actitud hacia las matemáticas, en consonancia con los objetivos del estudio.

La planificación se inició con la selección de actividades diseñadas para fomentar la formulación de preguntas, la creación de problemas matemáticos contextualizados y la resolución creativa de los mismos. Se eligieron situaciones cercanas y relevantes para los estudiantes, como celebraciones, ventas, emprendimientos, compras y mediciones, lo que facilitó la conexión con su entorno y sus experiencias cotidianas (ver Figuras I, II, III, IV y V). El objetivo principal de este enfoque fue despertar el interés y la curiosidad de los participantes, elementos clave para promover la creatividad y el pensamiento crítico (Stein, 2021).

En cuanto a los medios de observación, se prepararon cuadernillos de trabajo para los estudiantes, los cuales contenían las situaciones visuales propuestas y hojas específicas para cada fase del enfoque. Cada cuadernillo incluía una primera hoja con la imagen del problema, seguida de hojas dedicadas a las fases de exploración, creación y solución. Además, se proporcionaron materiales adicionales como hojas extra, reglas y calculadoras, según las necesidades de cada actividad, asegurando que los estudiantes contaran con los recursos necesarios para completar las tareas. El investigador y los docentes colaboradores también utilizaron un cuaderno de notas de campo.

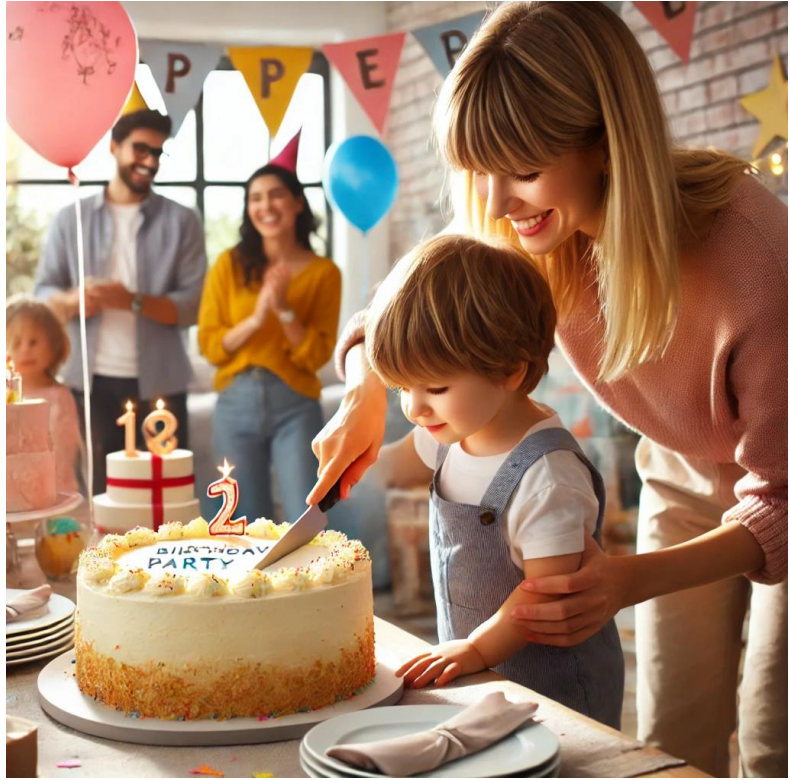

**Figura I**. Situación visual para 1° de secundaria

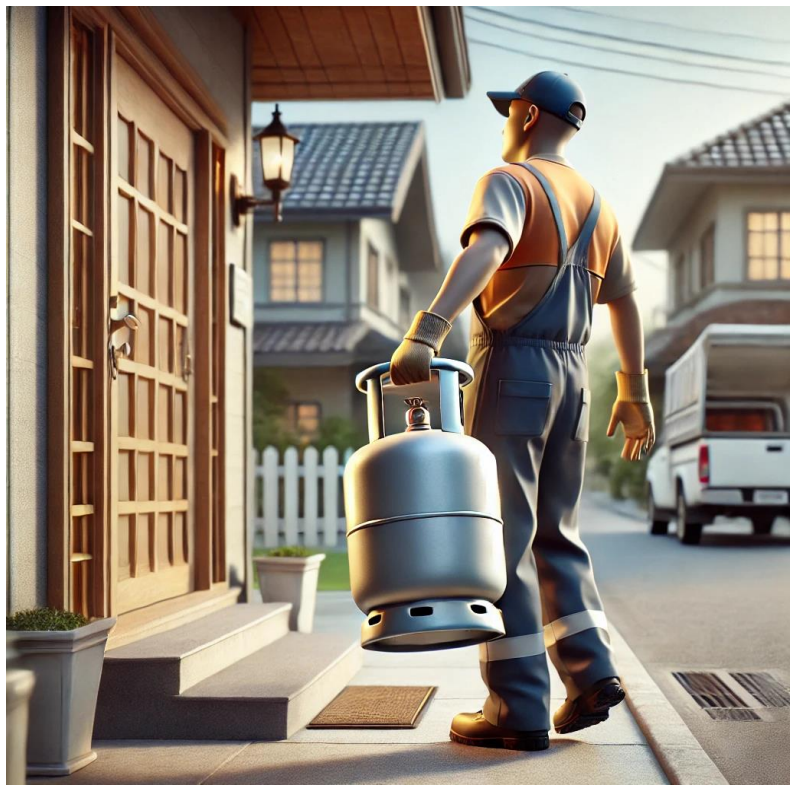

**Figura II**. Situación visual para 2° de secundaria

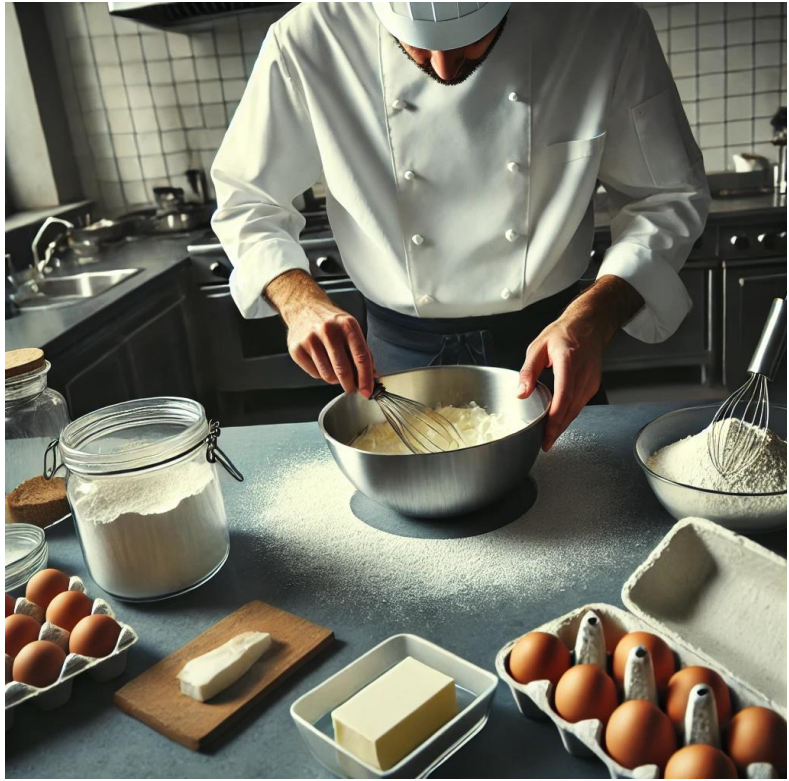

**Figura III**. Situación visual para 3° de secundaria

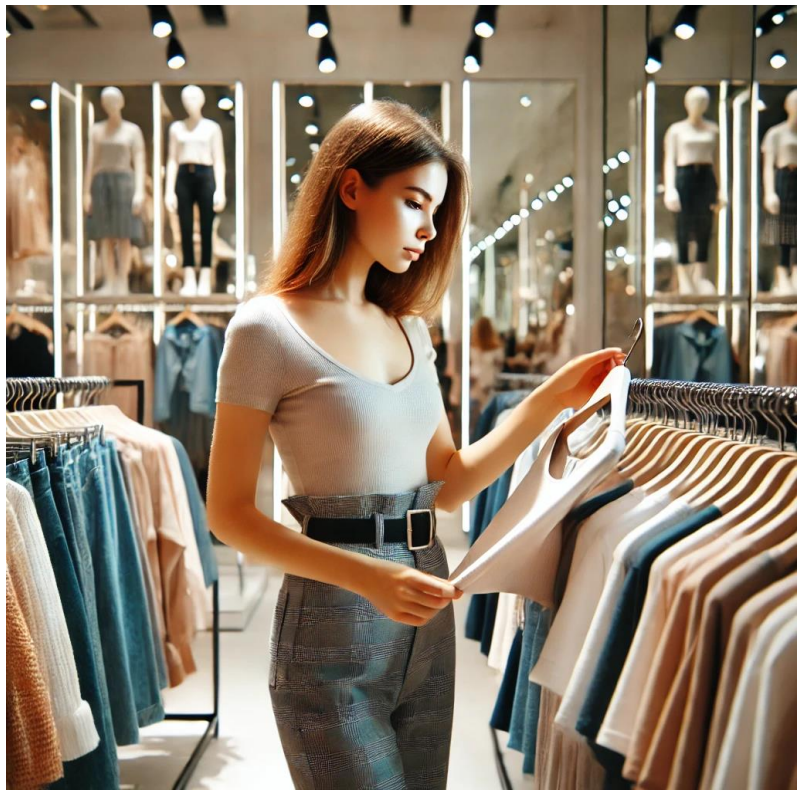

**Figura IV**. Situación visual para 4° de secundaria

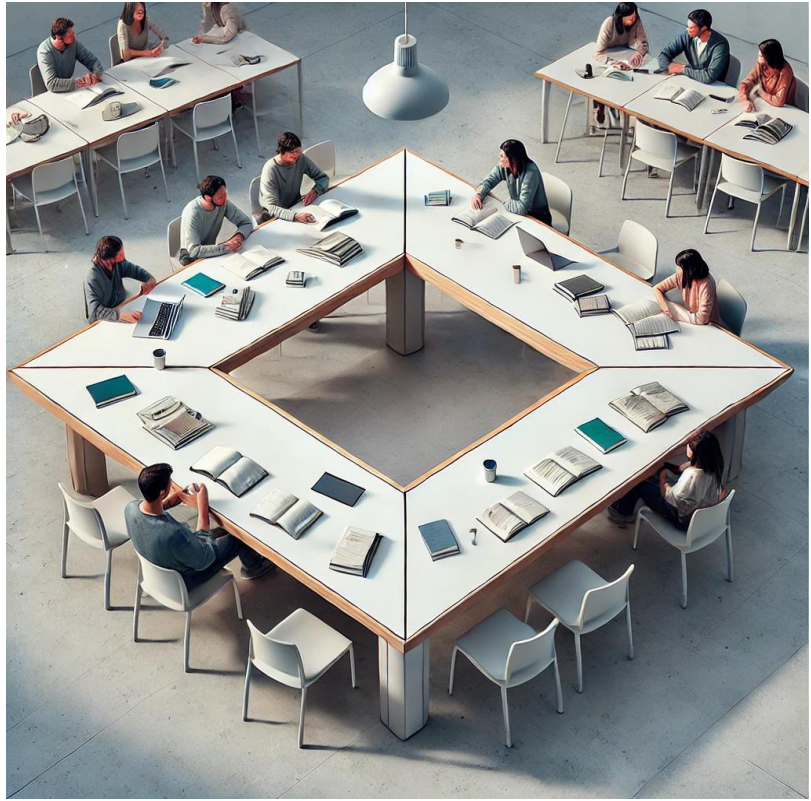

**Figura V**. Situación visual para 5° de secundaria

El desarrollo del concurso tuvo lugar en una única sesión de dos horas, durante la cual se llevaron a cabo las tres fases del enfoque adaptado. Esta estructura favoreció la continuidad del proceso de pensamiento y ayudó a los estudiantes a mantenerse enfocados en cada fase, permitiendo una inmersión completa en las actividades:

- Acto 1: Exploración. En esta fase, se presentó la situación problemática visual en el cuadernillo, pero sin ofrecer toda la información, con el fin de captar el interés de los estudiantes y fomentar su curiosidad. Los grupos generaron sus propias preguntas y seleccionaron la que les resultaba más interesante. Este proceso fue crucial para evaluar cómo los estudiantes aplicaban su creatividad y conocimientos previos, abordando el primer objetivo del estudio.
- Acto 2: Creación. Los grupos trabajaron de manera autónoma para construir un problema matemático coherente y contextualizado, basado en la pregunta seleccionada. Se hizo hincapié en la originalidad y la profundidad del problema, fomentando la flexibilidad cognitiva y el pensamiento divergente, habilidades clave en el proceso creativo según Guilford (1967). La colaboración entre los miembros del grupo facilitó el intercambio de ideas y estrategias, lo que permitió reinterpretar situaciones y ajustar enfoques cuando fue necesario, promoviendo así un enfoque más flexible y adaptativo en la resolución del problema.
- Acto 3: Solución. Cada grupo resolvió el problema que había creado, aplicando conceptos y razonamientos matemáticos pertinentes. La discusión en equipo facilitó que los estudiantes reflexionaran sobre sus métodos y resultados, y, en caso de inconsistencias, reformularan el problema. Este ejercicio promovió la reflexión metacognitiva, alentando a los estudiantes a evaluar y ajustar sus enfoques, lo que impulsó el desarrollo de habilidades de autorregulación y control cognitivo, tal como lo describe Flavell (1979).

El rol del investigador y los docentes fue esencial en este procedimiento. Actuaron como facilitadores del aprendizaje, guiando a los estudiantes sin intervenir directamente en sus decisiones, lo que les permitió tomar control de su propio proceso de aprendizaje, asumir la responsabilidad de sus soluciones y desarrollar autonomía y habilidades de autorregulación. Además, tanto los docentes como el investigador desempeñaron el rol de observadores durante el concurso, recopilando datos cualitativos mediante notas de campo y observación de la dinámica grupal. Esta observación fue clave para captar genuinamente las habilidades, estrategias y reflexiones de los estudiantes en relación con los objetivos específicos, sin interferir en su proceso de aprendizaje.

.

## 3.5. Técnicas de Recolección de Datos

Para obtener una comprensión integral del impacto de la adaptación del enfoque Matemáticas en Tres Actos, se emplearon dos técnicas principales de recolección de datos: la observación y la encuesta. La observación permitió captar las dinámicas grupales y los desempeños de los estudiantes durante las distintas fases del concurso, proporcionando datos cualitativos sobre sus interacciones y estrategias. Por otro lado, la encuesta, utilizando una escala de Likert, facilitó la obtención de datos cuantitativos relacionados con las percepciones de los estudiantes sobre su motivación, actitud y confianza en sus capacidades matemáticas, lo que complementó los hallazgos cualitativos y permitió una triangulación de datos más robusta.

**Observación de Datos**

En este estudio, la observación desempeñó un papel clave como técnica de recolección de datos, permitiendo registrar y analizar las dinámicas grupales y los desempeños de los estudiantes a lo largo de las fases del enfoque Matemáticas en Tres Actos. Esta técnica se aplicó de forma estructurada y deliberada, enfocándose en aspectos fundamentales del desempeño grupal, como la creatividad matemática, la resolución de problemas y las habilidades metacognitivas, tal como lo describe Angrosino (2012).

El proceso de observación se organizó en torno a cuatro elementos clave:

- El observador: el investigador principal y los docentes colaboradores.
- El ente observado: los desempeños colectivos de los grupos de estudiantes participantes.
- Los medios de observación: los cuadernillos de trabajo y los registros anecdóticos.
- Las circunstancias de la observación: el contexto del concurso de matemáticas.

El investigador principal, junto con los docentes colaboradores, asumió el rol de observadores participantes, integrándose al entorno del concurso para recopilar datos contextuales y detallados. Esta participación activa permitió observar las dinámicas grupales sin interferir en las decisiones o estrategias de los estudiantes, lo que garantizó la objetividad y autenticidad de las observaciones. Tal como señalan Supo y Zacarías (2020), la observación participante facilita la comprensión de las interacciones en su contexto natural, permitiendo que los datos sean recolectados sin alterar el desarrollo de las actividades.

Los desempeños colectivos de los grupos fueron evaluados a lo largo de las tres fases del enfoque adaptado: Exploración, Creación y Solución, constituyendo el objeto central del análisis. Estos desempeños proporcionaron información clave para entender el impacto de la metodología implementada en el aprendizaje de los estudiantes.

Para garantizar la sistematización y reproducibilidad de los hallazgos, se utilizaron como medios de observación los cuadernillos de trabajo y los registros anecdóticos. Los cuadernillos documentaron las preguntas formuladas, los problemas creados y las soluciones propuestas por los grupos, mientras que los registros anecdóticos capturaron observaciones cualitativas sobre las dinámicas grupales, las interacciones y las estrategias empleadas.

**Técnica Comunicacional: La Encuesta**

Para el componente cuantitativo del estudio, se empleó la encuesta como técnica principal de recolección de datos. El instrumento utilizado fue una escala de Likert, diseñada específicamente para evaluar las percepciones de los participantes en relación con tres dimensiones clave: motivación, actitud y confianza en sus capacidades matemáticas, tras su participación en el concurso basado en la adaptación del enfoque Matemáticas en Tres Actos. La justificación de esta elección se basa en la capacidad de las encuestas para recolectar datos de manera estandarizada y eficiente, permitiendo obtener información directa de los participantes sobre sus experiencias y percepciones. Como señalan Cohen, Manion y Morrison (2018), las encuestas son herramientas ideales para explorar aspectos subjetivos de la experiencia individual, proporcionando una base sólida para el análisis cuantitativo y la interpretación de tendencias generales.

## 3.6. Instrumentos y Medios de Observación

La recolección de datos en este estudio se basó en una combinación de observaciones directas, análisis de cuadernillos de trabajo y la aplicación de un cuestionario, con el objetivo de evaluar en profundidad los aspectos relacionados con los objetivos específicos del estudio.

**Medios de Observación: Cuadernillos de Trabajo y Registros Anecdóticos**

Los cuadernillos de trabajo utilizados durante el concurso fueron una fuente clave de información. Estos documentos contenían las producciones escritas de los estudiantes, desde las preguntas formuladas en el Acto 1 hasta las soluciones finales en el Acto 3. El análisis de estos cuadernillos permitió examinar en profundidad los problemas creados por los estudiantes, evaluando la coherencia en sus formulaciones y la originalidad en sus respuestas.

Adicionalmente, se utilizaron registros anecdóticos como medio complementario de observación. Estos registros capturaron detalles cualitativos sobre las dinámicas grupales, las interacciones y las estrategias empleadas por los estudiantes durante las fases del enfoque adaptado. Los registros anecdóticos proporcionaron una visión más rica y contextualizada de los procesos de aprendizaje.

**Instrumento: Escala de Likert**

Para el componente cuantitativo del estudio, se empleó una escala de Likert como instrumento principal para la recolección de datos. Esta escala fue diseñada específicamente con el fin de medir las percepciones de los estudiantes en tres dimensiones clave: motivación, actitud y confianza en sus habilidades matemáticas, tras su participación en el concurso.

La escala incluyó afirmaciones relacionadas con las tres dimensiones evaluadas. Cada afirmación fue cuidadosamente formulada para reflejar con precisión las percepciones de los estudiantes, y las respuestas se evaluaron en una escala de cinco puntos, que iba desde "Totalmente en desacuerdo" hasta "Totalmente de acuerdo".

## 3.7 Análisis de Datos

Este apartado describe el proceso de análisis cualitativo y cuantitativo empleado en el estudio para evaluar el impacto de la adaptación del enfoque Matemáticas en Tres Actos en las dimensiones clave: creatividad matemática, resolución de problemas y habilidades metacognitivas. Además, se analiza la percepción de los estudiantes hacia el aprendizaje matemático tras la implementación del enfoque.

**Definición de Indicadores y Codificación**

Para evaluar los objetivos específicos del estudio, se definieron indicadores precisos para cada dimensión, lo que permitió una evaluación detallada de las competencias y habilidades desarrolladas por los estudiantes. Cada indicador representó una competencia esencial y fue codificado para facilitar el análisis de los datos tanto cualitativos como cuantitativos.

a) **Dimensión: Creatividad Matemática**
- DEFFP: Diversidad de Enfoques y Fluidez en la Formulación de Preguntas
- ICC: Interpretación Creativa del Contexto
- IEP: Innovación en el Enfoque del Problema
- CP: Contextualización del Problema
- CSD: Calidad y Suficiencia de los Datos
- PVP: Pertinencia y Variedad de las Preguntas
- EEM: Existencia de un Entorno Matemático

b) **Dimensión: Resolución de Problemas**
- CCP: Comunica su Comprensión del Problema
- TPLM: Traduce Problemas a Lenguaje Matemático
- PEP: Pertinencia de Estrategias y Procedimientos
- JPR: Justificación de Procedimientos y Resultados

c) **Dimensión: Habilidades Metacognitivas**
- ATD: Autonomía en la Toma de Decisiones
- FEP: Flexibilidad en la Estructuración del Problema
- ICE: Identificación y Corrección de Errores
- AAE: Autorregulación y Ajuste de Estrategias

Estos indicadores fueron codificados para facilitar su evaluación y permitir el análisis tanto cualitativo como cuantitativo de los datos recopilados durante el estudio.

**Niveles de Desempeño por Indicador**

Para cada indicador, se establecieron cuatro niveles de desempeño: Destacado, Satisfactorio, En Proceso y En Inicio. Estos niveles permitieron una evaluación estandarizada y detallada del desarrollo de competencias en cada área, facilitando la comparación y análisis de los resultados entre los grupos y a lo largo de las diferentes fases del estudio.

**Alineación de Indicadores y Objetivos Específicos**

Los indicadores se agruparon en tres dimensiones principales, alineadas con los objetivos específicos del estudio, lo que estructuró el análisis de manera coherente. Además, se consideró una cuarta dimensión para explorar las percepciones de los estudiantes, evaluadas mediante la escala de Likert utilizada en la encuesta. Esta organización permitió una evaluación integral del impacto de la intervención en diferentes áreas del aprendizaje matemático.

**Procedimiento de Análisis de Datos**

El análisis de datos se llevó a cabo mediante un enfoque mixto, combinando técnicas cualitativas y cuantitativas, lo que permitió una evaluación completa y profunda del impacto de la adaptación del enfoque en las diferentes dimensiones del aprendizaje.

a) **Análisis cualitativo**: Este análisis se centró en las producciones de los estudiantes, que incluyeron los cuadernillos de trabajo y las observaciones cualitativas obtenidas a través de registros anecdóticos. El proceso comenzó con la recopilación de datos y una codificación inicial, en la que se revisaron las producciones para identificar evidencias relacionadas con los indicadores definidos, asignándose los códigos correspondientes, de acuerdo con Miles, Huberman y Saldaña (2014).

Posteriormente, se evaluaron los indicadores para cada grupo, asignando un nivel de desempeño según los criterios establecidos. Este proceso facilitó la comparación de los desempeños entre los diferentes grupos y grados, permitiendo la identificación de patrones, fortalezas y áreas de mejora.

b) **Análisis cuantitativo**: En el análisis cuantitativo, se enfocó en las percepciones de los estudiantes, las cuales fueron analizadas mediante los datos obtenidos de la escala de Likert aplicada al finalizar la intervención. Tras la recolección de las respuestas, se realizaron análisis estadísticos descriptivos, calculando frecuencias y promedios para cada ítem de la encuesta. Además, se identificaron tendencias y patrones en las percepciones de los estudiantes hacia el enfoque y su aprendizaje matemático.

**Interpretación y Validación del Análisis**

La combinación de los análisis cualitativo y cuantitativo permitió una evaluación integral del impacto del enfoque adaptado en las diferentes dimensiones de aprendizaje. Esta triangulación de datos fortaleció la confiabilidad del estudio, proporcionando una visión más rica y sólida de los resultados obtenidos. Se contrastaron los hallazgos de ambos enfoques (cualitativo y cuantitativo) para validar los resultados y examinar la coherencia entre las observaciones de desempeño y las percepciones reportadas por los estudiantes.

**Consideraciones Éticas**

Se garantizó la confidencialidad y el anonimato de los participantes a lo largo de todo el proceso de recolección y análisis de datos. La participación en el estudio fue voluntaria, y se obtuvo el consentimiento informado de los estudiantes y, cuando fue necesario, de sus tutores legales. Además, se respetaron los principios éticos de la investigación educativa, asegurando que los datos recolectados fueran utilizados exclusivamente con fines académicos y para la mejora pedagógica.

# 4. Resultados

Se presentan los resultados del estudio sobre la aplicación del enfoque adaptado de Matemáticas en Tres Actos, organizado en cuatro dimensiones: creatividad matemática, resolución de problemas, habilidades metacognitivas y percepciones y actitudes hacia el aprendizaje matemático. Se analiza cómo los estudiantes de diferentes grados aplicaron la estrategia en las tres fases del enfoque, evaluando su impacto en el desempeño académico y las percepciones de los estudiantes, con ejemplos representativos y datos respaldados por tablas y gráficos.

## 4.1 Resultados dimensión Creatividad Matemática

En la Tabla I se presentan los resultados obtenidos por los grupos participantes en la dimensión de Creatividad Matemática. Los datos muestran el desempeño en siete indicadores clave: Diversidad de Enfoques y Fluidez en la Formulación de Preguntas (DEFFP), Innovación en la Creación de Conceptos (ICC), Innovación en el Enfoque del Problema (IEP), Contextualización del Problema (CP), Calidad y Suficiencia de los Datos (CSD), Pertinencia y Variedad de las Preguntas (PVP) y Existencia de un Entorno Matemático (EEM). Estos resultados, clasificados en las escalas de Destacado, Satisfactorio, En Proceso y En Inicio, permiten observar con detalle las fortalezas y áreas de mejora en cada indicador.

**Tabla I**. Resultados dimensión Creatividad Matemática

| Grado | Grupo | DEFFP | ICC | IEP | CP | CSD | PVP | EEM |
|---|---|---|---|---|---|---|---|---|
| Primero | P-01 | Satisfactorio | Satisfactorio | En Proceso | Destacado | Destacado | Satisfactorio | Destacado |
| | P-02 | Destacado | Destacado | Destacado | Destacado | Satisfactorio | Satisfactorio | Destacado |
| | P-03 | En Proceso | Satisfactorio | En Proceso | Satisfactorio | En Proceso | En Proceso | Satisfactorio |
| | P-04 | Satisfactorio | En Proceso | En Proceso | En Proceso | En Inicio | En Proceso | En Proceso |
| Segundo | S-01 | Satisfactorio | Destacado | Destacado | Satisfactorio | Destacado | Satisfactorio | Destacado |
| | S-02 | En Proceso | Satisfactorio | Satisfactorio | Destacado | Satisfactorio | En Proceso | En Inicio |
| | S-03 | Destacado | En Proceso | En Proceso | Satisfactorio | En Proceso | Satisfactorio | Destacado |
| | S-04 | En Proceso | En Proceso | En Proceso | Satisfactorio | En Proceso | Satisfactorio | En Proceso |
| Tercero | T-01 | Satisfactorio | En Proceso | En Proceso | Satisfactorio | En Proceso | Satisfactorio | Satisfactorio |
| | T-02 | En Proceso | En Proceso | En Proceso | En Proceso | En Proceso | En Proceso | Satisfactorio |
| | T-03 | En proceso | En proceso | En inicio | En proceso | Satisfactorio | En proceso | Satisfactorio |
| Cuarto | C-01 | Destacado | Satisfactorio | En Proceso | Destacado | En Proceso | En Proceso | Destacado |
| | C-02 | Satisfactorio | Destacado | Destacado | Satisfactorio | En Proceso | Satisfactorio | Destacado |
| | C-03 | Satisfactorio | Satisfactorio | En proceso | En proceso | En proceso | En proceso | Destacado |
| | C-04 | En Proceso | Satisfactorio | Satisfactorio | En Proceso | En Proceso | En Proceso | Destacado |
| | C-05 | Satisfactorio | Satisfactorio | En Proceso | Satisfactorio | En Inicio | En Proceso | Satisfactorio |
| Quinto | Q-01 | En Proceso | Destacado | Destacado | Satisfactorio | En Proceso | En Proceso | Destacado |
| | Q-02 | En Proceso | Satisfactorio | En Proceso | Satisfactorio | En Proceso | En Proceso | Satisfactorio |

En términos de diversidad de enfoques y fluidez en la formulación de preguntas, la mayoría de los estudiantes, especialmente en los grados superiores, logró un desempeño destacado al formular preguntas complejas y variadas. Sin embargo, algunos grupos de grados inferiores se limitaron a enfoques más concretos, lo que sugiere la necesidad de fomentar la creatividad desde las primeras fases del proceso.

En cuanto a la interpretación creativa del contexto (ICC) e innovación en el enfoque del problema (IEP), los estudiantes demostraron su capacidad para reimaginar contextos y adoptar enfoques originales. No obstante, algunos grupos necesitaron mayor apoyo para estructurar problemas innovadores de manera efectiva.

Respecto a la coherencia del problema creado, los estudiantes lograron contextualizar bien los problemas, aunque algunos grupos enfrentaron dificultades para asegurar la congruencia de los datos y formular preguntas variadas. A pesar de esto, la mayoría mostró una sólida capacidad para aplicar conceptos matemáticos relevantes, situando los problemas en un entorno matemático adecuado y reflejando su comprensión de los conceptos.

Para complementar estos resultados, el gráfico mostrado en la Figura VI, resume los promedios de desempeño en los siete indicadores evaluados, proporcionando una visión global del rendimiento de los grupos. Los valores promedio, clasificados en una escala de En Inicio (1) a Destacado (4), permiten observar la distribución de los desempeños en cada indicador clave.

El análisis del gráfico destaca tanto las fortalezas como las áreas de mejora en la dimensión de creatividad matemática. Por ejemplo, el indicador EEM alcanzó el promedio más alto, lo que sugiere que los grupos fueron efectivos al integrar entornos matemáticos coherentes en sus propuestas. En contraste, el indicador PVP mostró el promedio más bajo, lo que señala dificultades en la formulación de preguntas variadas y pertinentes. Este análisis no solo complementa los hallazgos de la tabla, sino que también refuerza los objetivos del estudio, proporcionando una visión clara de cómo se manifiestan las distintas facetas de la creatividad matemática en el contexto del concurso.

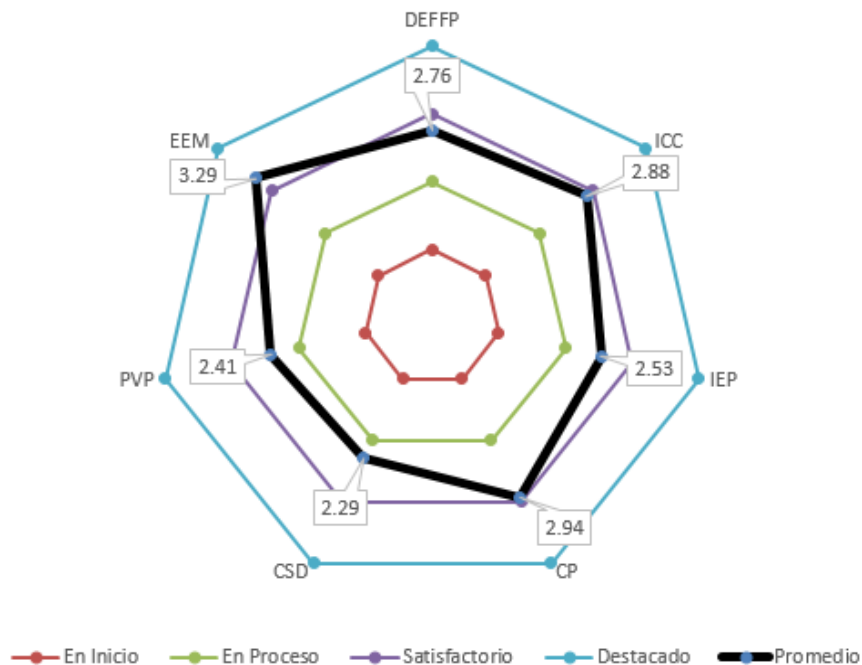

**Figura VI**. Resultados dimensión Creatividad Matemática

A continuación, se presentan análisis detallados de grupos representativos de cada grado, ilustrando cómo aplicaron la creatividad matemática en relación con los indicadores clave:

- Primer Grado (Grupo P-02): Este grupo demostró una gran creatividad al formular preguntas originales sobre temas como el MCM, áreas y teoría de conjuntos. Introdujeron un contexto creativo, aplicando la teoría de conjuntos de manera innovadora, alcanzando un enfoque destacado en los indicadores de fluidez (DEFFP), interpretación creativa (ICC) e innovación en el enfoque del problema (IEP). A pesar de una inconsistencia inicial en los datos, lograron corregir el error con precisión.
- Segundo Grado (Grupo S-01): Este grupo formuló preguntas claras relacionadas con el movimiento rectilíneo uniforme y el MCM. Introdujeron un enfoque novedoso al aplicar el MCM a la periodicidad de los pedidos de gas, alcanzando niveles destacados en varios indicadores.
- Tercer Grado (Grupo T-01): Los estudiantes de este grupo formularon preguntas sobre cálculos de precios y ganancias, pero representaron precios mediante polinomios. Aunque exploraron conceptos avanzados, la abstracción redujo la aplicabilidad práctica del problema, lo que sugiere la necesidad de equilibrar la innovación con la relevancia práctica.
- Cuarto Grado (Grupo C-02): Este grupo mostró gran creatividad al introducir un concepto único de "manía adictiva" de compras, aplicando progresiones cuadráticas para modelar el patrón de compras. Sin embargo, algunos detalles fueron omitidos en la formulación inicial.
- Quinto Grado (Grupo Q-01): Introdujeron una cucaracha personificada en una disposición en forma de "U", aplicando teoría de grafos para determinar rutas. Aunque demostraron creatividad, la falta de precisión en los datos limitó el desarrollo del análisis matemático.

En resumen, los resultados en la dimensión de Creatividad Matemática están estrechamente ligados al primer objetivo del estudio: evaluar la capacidad de los estudiantes para generar y formular problemas matemáticos creativos. El enfoque implementado fomentó la creatividad, permitiendo a los estudiantes reinterpretar contextos y aplicar conceptos matemáticos de manera innovadora. Grupos como P-02, S-01, C-02 y Q-01 se destacaron por su capacidad para innovar, lo que demuestra que la metodología promueve un ambiente favorable para el pensamiento creativo y la exploración más allá de enfoques tradicionales. Además, su habilidad para conectar conceptos teóricos con situaciones prácticas subraya la efectividad del método para desarrollar una comprensión profunda y aplicada de las matemáticas, crucial para preparar a los estudiantes a resolver problemas reales.

## 4.2 Resultados Dimensión Resolución de Problemas

La Tabla II Resultados en la dimensión Resolución de Problemas presenta el desempeño de los grupos participantes en los cuatro indicadores evaluados: Comunica su Comprensión del Problema (CCP), Traduce Problemas a Lenguaje Matemático (TPLM), Pertinencia de Estrategias y Procedimientos (PEP) y Justificación de Procedimientos y Resultados (JPR). Los resultados están organizados por grado y grupo, y clasificados en las escalas de Destacado, Satisfactorio, En Proceso y En Inicio, lo que permite identificar fortalezas y áreas de mejora en las distintas fases de la resolución de problemas.

**Tabla II**. Resultados en la dimensión Resolución de Problemas

| Grado | Grupo | CCP | TPLM | PEP | JPR |
|---|---|---|---|---|---|
| Primero | P-01 | Destacado | Destacado | Destacado | Satisfactorio |
| | P-02 | Destacado | Destacado | Destacado | Destacado |
| | P-03 | Satisfactorio | Satisfactorio | Destacado | Satisfactorio |
| | P-04 | En Proceso | En proceso | En proceso | En proceso |
| Segundo | S-01 | Destacado | Destacado | Destacado | Destacado |
| | S-02 | Satisfactorio | Satisfactorio | Satisfactorio | Satisfactorio |
| | S-03 | Satisfactorio | Destacado | Destacado | Destacado |
| | S-04 | Satisfactorio | Satisfactorio | Satisfactorio | En proceso |
| Tercero | T-01 | Satisfactorio | Satisfactorio | Satisfactorio | Satisfactorio |
| | T-02 | Satisfactorio | Satisfactorio | Satisfactorio | Satisfactorio |
| | T-03 | En proceso | Satisfactorio | Satisfactorio | Satisfactorio |
| Cuarto | C-01 | Satisfactorio | Destacado | Destacado | Satisfactorio |
| | C-02 | Destacado | Satisfactorio | Destacado | Satisfactorio |
| | C-03 | Satisfactorio | Satisfactorio | Satisfactorio | Satisfactorio |
| | C-04 | Satisfactorio | Destacado | Satisfactorio | En proceso |
| | C-05 | En proceso | En proceso | En Inicio | En Inicio |
| Quinto | Q-01 | Destacado | Satisfactorio | Destacado | Satisfactorio |
| | Q-02 | Satisfactorio | Satisfactorio | En proceso | En proceso |

Los estudiantes, en su mayoría, destacaron en CCP y TPLM, lo que indica que fueron capaces de comprender y expresar los problemas matemáticos con claridad. Además, en el indicador PEP, los estudiantes demostraron una sólida capacidad para aplicar estrategias y procedimientos matemáticos adecuados, alcanzando niveles satisfactorios o destacados. Sin embargo, se identificaron áreas de mejora en JPR, donde la justificación de procedimientos y resultados no siempre fue tan detallada o rigurosa como se esperaba. Algunos grupos ofrecieron explicaciones superficiales, lo que sugiere que hay espacio para fortalecer la habilidad de los estudiantes para argumentar matemáticamente y profundizar en la justificación de sus soluciones.

Complementariamente, el gráfico radar (Figura VII) resume los resultados promedio obtenidos en los cuatro indicadores evaluados, proporcionando una visión global del desempeño colectivo en esta dimensión. Los valores promedio, clasificados de En Inicio (1) a Destacado (4), permiten observar cómo se distribuyen los desempeños en cada indicador clave. El análisis del gráfico revela que los indicadores TPLM y PEP obtuvieron los promedios más altos (3.22), lo que sugiere que los estudiantes fueron efectivos en traducir problemas a lenguaje matemático y en la aplicación de estrategias pertinentes. En contraste, el indicador JPR, con un promedio de 2.83, muestra un desempeño más bajo, lo que destaca la necesidad de fortalecer la capacidad de los estudiantes para justificar de manera más sólida sus procedimientos y resultados.

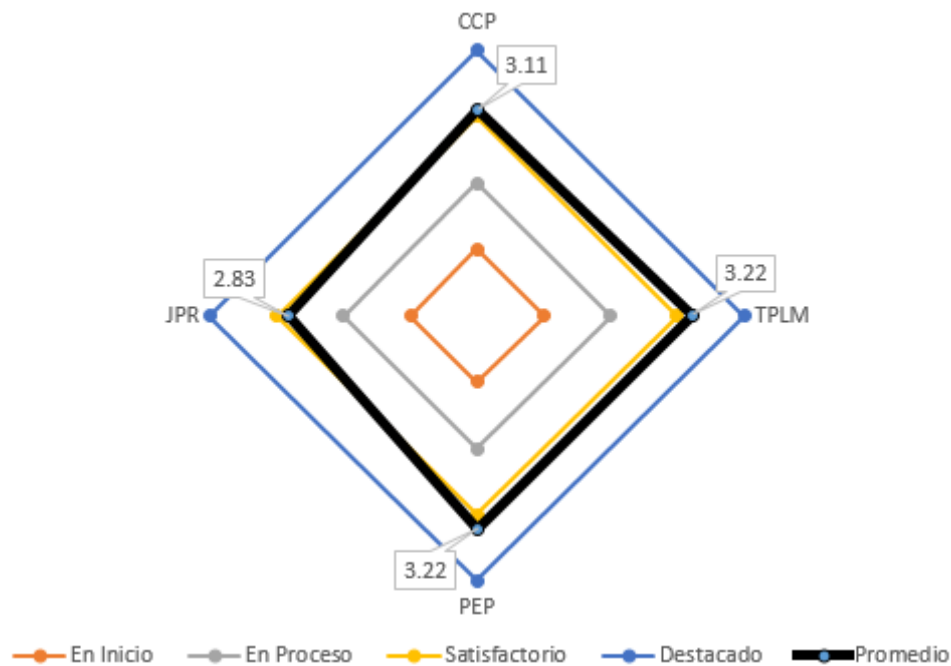

**Figura VII**. Resultados en la dimensión Resolución de Problem*as*

A continuación, se ofrecen análisis detallados de grupos representativos de cada grado, destacando cómo aplicaron la creatividad matemática en relación con los indicadores clave:

- Primer Grado (Grupo P-02): Este grupo mostró un desempeño destacado en todos los indicadores, especialmente en CCP y TPLM, donde demostraron una comprensión profunda del problema y lo tradujeron de manera efectiva al lenguaje matemático. En PEP, aplicaron estrategias eficientes, pero en JPR pudieron haber profundizado más en la explicación de sus decisiones estratégicas.
- Segundo Grado (Grupo S-01): El grupo S-01 destacó en todos los indicadores. Su capacidad para comunicar la comprensión del problema y traducirlo al lenguaje matemático fue excepcional. También aplicaron estrategias adecuadas y justificaron sus procedimientos de manera completa, mostrando un alto nivel de comprensión matemática.
- Tercer Grado (Grupo T-01): El grupo T-01 alcanzó un desempeño satisfactorio en CCP y TPLM, pero las justificaciones en JPR fueron más limitadas. Aunque aplicaron estrategias matemáticas adecuadas en PEP, su razonamiento detrás de los procedimientos no fue lo suficientemente detallado, lo que resalta la necesidad de fortalecer la argumentación matemática.
- Cuarto Grado (Grupo C-02): Este grupo presentó un desempeño destacado, especialmente en CCP, donde desglosaron el problema de manera clara y detallada. En TPLM, tradujeron el problema matemáticamente utilizando expresiones algebraicas para representar precios, descuentos y ganancias, facilitando su resolución. Las justificaciones en JPR fueron detalladas, lo que demuestra una comprensión profunda del proceso matemático involucrado.

- Quinto Grado (Grupo Q-01): El grupo mostró un desempeño satisfactorio en la resolución de problemas. Aunque comprendieron el problema y lo representaron utilizando teoría de grafos, su análisis en JPR fue limitado, lo que sugiere que se debe fomentar una mayor profundización en la justificación matemática.

Los resultados en la dimensión de Resolución de Problemas están directamente vinculados con el segundo objetivo específico del estudio: analizar el desempeño de los estudiantes en la resolución de problemas matemáticos. El enfoque implementado contribuyó al desarrollo de habilidades analíticas, permitiendo a los estudiantes comprender problemas complejos y traducirlos al lenguaje matemático. Grupos como P-02, S-01 y C-02 demostraron capacidad para descomponer problemas y representarlos de manera estructurada.

## 4.3 Resultados en la Dimensión Habilidades Metacognitivas

La Tabla III presenta una visión detallada del desempeño de los grupos participantes en la dimensión Habilidades Metacognitivas, desglosando los resultados por grado y grupo en los indicadores evaluados: Autonomía en la Toma de Decisiones (ATD), Flexibilidad en la Estructuración (FEP), Identificación y Corrección de Errores (ICE) y Autorregulación y Ajuste de Estrategias (AAE). Esta información permite identificar claramente los patrones de desempeño, destacando las fortalezas y las áreas de mejora específicas para cada grupo.

**Tabla III**. Resultados en la dimensión Habilidades Metacognitivas

| Grado | Grupo | ATD | FEP | ICE | AAE |
|---|---|---|---|---|---|
| Primero | P-01 | Satisfactorio | Satisfactorio | Satisfactorio | Satisfactorio |
| | P-02 | Destacado | Destacado | Destacado | Destacado |
| | P-03 | Satisfactorio | En proceso | Satisfactorio | Satisfactotio |
| | P-04 | En Proceso | En Inicio | En Inicio | En Inicio |
| Segundo | S-01 | Destacado | Satisfactorio | Destacado | Destacado |
| | S-02 | Satisfactorio | Satisfactorio | Satisfactorio | Satisfactorio |
| | S-03 | Satisfactorio | Satisfactorio | Satisfactorio | Satisfactorio |
| | S-04 | En Proceso | En Proceso | En Proceso | En proceso |
| Tercero | T-01 | Satisfactorio | En proceso | Destacado | Satisfactorio |
| | T-02 | Satisfactorio | En proceso | En proceso | En proceso |
| | T-03 | En proceso | En Inicio | En Inicio | En Inicio |
| Cuarto | C-01 | Destacado | Satisfactorio | Satisfactorio | Satisfactorio |
| | C-02 | Destacado | Satisfactorio | Destacado | Satisfactorio |
| | C-03 | Destacado | En proceso | En proceso | Satisfactorio |
| | C-04 | Destacado | En proceso | En proceso | Satisfactorio |
| | C-05 | En proceso | En inicio | En inicio | En inicio |
| Quinto | Q-01 | Satisfactorio | Satisfactorio | En proceso | Satisfactorio |
| | Q-02 | Satisfactorio | En proceso | En proceso | En proceso |

La mayoría mostró un desempeño satisfactorio o destacado en Autonomía en la Toma de Decisiones (ATD), al tomar decisiones de forma independiente al formular y resolver problemas, lo que refleja confianza en sus capacidades y disposición para asumir responsabilidad en su aprendizaje. En Flexibilidad en la Estructuración del Problema (FEP), algunos grupos demostraron la capacidad de reestructurar problemas al identificar inconsistencias, mientras que otros se mantuvieron rígidos, lo que sugiere la necesidad de fomentar mayor flexibilidad cognitiva. En cuanto a Identificación y Corrección de Errores (ICE), se observó que algunos estudiantes detectaron y corrigieron inconsistencias, pero otros pasaron por alto errores que afectaron la precisión y coherencia de sus soluciones, lo que resalta la oportunidad de mejorar la revisión crítica y autoverificación. Finalmente, en Autorregulación y Ajuste de Estrategias (AAE), los estudiantes ajustaron sus estrategias cuando enfrentaron obstáculos, mostrando autorregulación, pero estos ajustes fueron a veces superficiales, lo que indica que se requiere una reflexión más profunda sobre las dificultades encontradas.

El gráfico radar (Figura VIII), sintetiza los resultados promedio en los cuatro indicadores evaluados, proporcionando una visión general de la dimensión de habilidades metacognitivas. Los datos muestran que Autonomía en la Toma de Decisiones (ATD) obtuvo el mejor desempeño promedio con un valor de 3.11, lo que refleja que los estudiantes lograron tomar decisiones autónomas de manera efectiva durante el concurso. En cambio, Flexibilidad en la Estructuración del Problema (FEP) registró el promedio más bajo, con un valor de 2.33, lo que destaca las dificultades en la reestructuración flexible de estrategias. Los indicadores Identificación y Corrección de Errores (ICE) y Autorregulación y Ajuste de Estrategias (AAE) mostraron promedios intermedios de 2.56 y 2.59, respectivamente, lo que sugiere un desempeño moderado en la capacidad de corregir errores y ajustar estrategias. Este análisis gráfico permite visualizar claramente las fortalezas y áreas críticas de mejora en las habilidades metacognitivas, alineándose con los objetivos del estudio y proporcionando información clave para futuras implementaciones educativas.

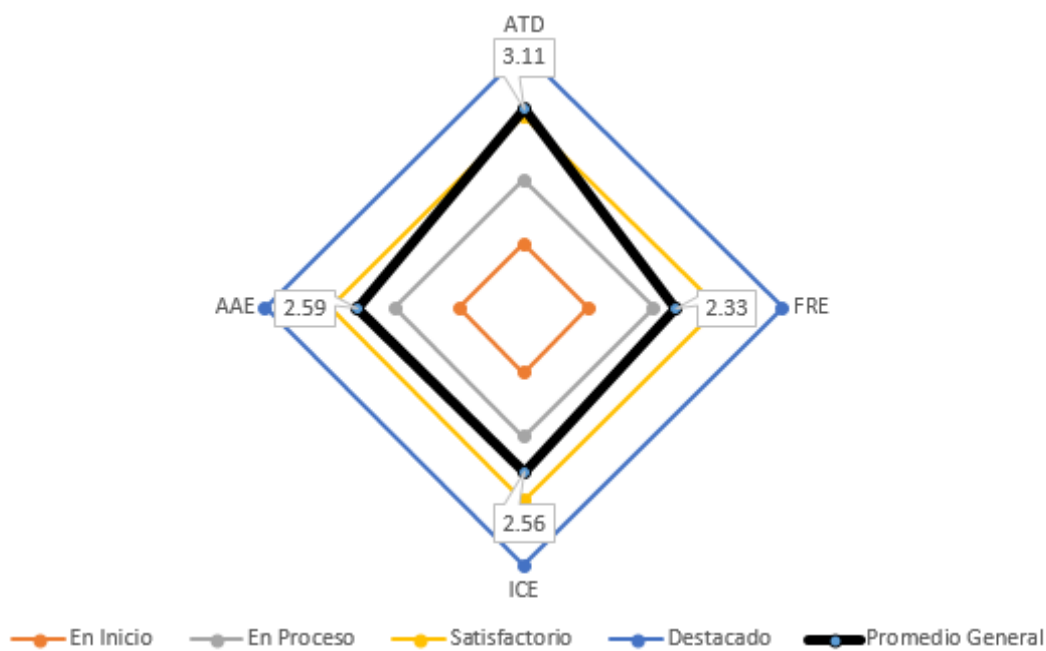

**Figura VIII**. Resultados en la dimensión Habilidades Metacognitivas

A continuación, se presentan ejemplos representativos de grupos de cada grado, ilustrando su desempeño en los indicadores clave de habilidades metacognitivas:

- Primer Grado (Grupo P-02): El Grupo demostró un desempeño destacado en Autonomía en la Toma de Decisiones (ATD), tomando la iniciativa en la formulación y resolución del problema. En Flexibilidad en la Estructuración del Problema (FEP), mostraron flexibilidad al reestructurar el problema al identificar inconsistencias en los datos. No obstante, en Identificación y Corrección de Errores (ICE), aunque corrigieron el error en los datos, podrían mejorar en la detección temprana de inconsistencias. En Autorregulación y Ajuste de Estrategias (AAE), ajustaron sus estrategias de manera efectiva para llegar a la solución correcta.
- Segundo Grado (Grupo S-01): El Grupo mostró un desempeño destacado en ATD, desarrollando y resolviendo el problema de manera independiente. En FEP, mantuvieron una estructura coherente sin necesidad de reestructuración significativa. En ICE, no se observaron errores significativos, lo que refleja precisión en su trabajo. En AAE, planificaron y ejecutaron sus estrategias de manera efectiva sin requerir ajustes mayores.
- Tercer Grado (Grupo T-01): El Grupo tuvo un desempeño satisfactorio en ATD, tomando decisiones en la resolución del problema. En FEP, se mostraron rígidos, sin adaptar el planteamiento a pesar de las dificultades, lo que indica una oportunidad para mejorar en flexibilidad. En ICE, pasaron por alto inconsistencias que afectaron la aplicabilidad práctica del problema. En AAE, realizaron ajustes superficiales en sus estrategias sin una reflexión profunda sobre las dificultades enfrentadas.
- Cuarto Grado (Grupo C-02): El Grupo demostró un desempeño satisfactorio en habilidades metacognitivas. En ATD, mostraron independencia al tomar decisiones para formular y resolver el problema, seleccionando estrategias apropiadas. En FEP, realizaron algunos ajustes menores durante el proceso de resolución, indicando flexibilidad en la estructuración del problema. En ICE, inicialmente omitieron detalles importantes en los datos, pero lograron identificar y corregir estas omisiones durante la resolución, mostrando capacidad para rectificar errores. En AAE, ajustaron sus estrategias frente a las dificultades, demostrando autorregulación y adaptación, aunque podrían beneficiarse de una reflexión más profunda sobre los ajustes realizados.
- Quinto Grado (Grupo Q-01): El Grupo presentó un desempeño satisfactorio en habilidades metacognitivas. En ATD, tomaron decisiones autónomas al desarrollar y resolver el problema, incluyendo ajustes cuando fue necesario. En FEP, mostraron flexibilidad al realizar adaptaciones parciales durante el proceso, aunque podrían haber sido más abiertos a reestructurar su enfoque para mejorar la precisión y claridad del problema. En ICE, pasaron por alto la falta de datos precisos en el enunciado, lo que afectó la exactitud de su solución, indicando la necesidad de mejorar en la identificación y corrección de errores. En AAE, ajustaron sus estrategias durante la resolución, mostrando capacidad de autorregulación, aunque una reflexión más profunda sobre sus procesos podría haber mejorado su desempeño.

Los resultados en la dimensión de Habilidades Metacognitivas están estrechamente vinculados con el tercer objetivo específico del estudio: examinar el desarrollo de habilidades metacognitivas en el proceso de creación y resolución de problemas. Los hallazgos muestran fortalezas y áreas de mejora, destacando que el enfoque metodológico promovió el desarrollo de la autonomía y autorregulación, observándose un alto grado de independencia en grupos como P-02, S-01 y C-02. Sin embargo, se identificó una necesidad de fortalecer la flexibilidad cognitiva, especialmente en los grupos de tercer y quinto grado, quienes mostraron resistencia a reestructurar problemas y explorar alternativas.

## 4.4 Percepciones y Actitudes hacia el Aprendizaje Matemático

Se analizaron las percepciones y actitudes de los estudiantes hacia las matemáticas tras su participación en el concurso, donde se implementó el enfoque adaptado de Matemáticas en Tres Actos. Este análisis se basó en los resultados de un cuestionario administrado al final de la intervención, utilizando una escala de Likert para evaluar tres aspectos clave: motivación, actitud hacia las matemáticas y confianza en su capacidad matemática. Los resultados proporcionan una visión global del impacto del enfoque en la experiencia educativa de los participantes.

La evaluación de estos tres indicadores se realizó considerando tanto el puntaje promedio por pregunta como los niveles alcanzados en cada uno. Utilizando una escala de Likert de 1 a 5, se midieron diversos aspectos relacionados con el aprendizaje matemático, lo que permitió obtener una visión completa de la percepción de los estudiantes.

**Resultados relacionados con la Motivación**

Los resultados muestran una tendencia general positiva, con un hallazgo destacado en la Pregunta 2, sobre el interés por continuar aprendiendo matemáticas, que obtuvo el puntaje promedio más alto (3.58). Esto indica que las actividades del enfoque generaron un interés renovado en los estudiantes. Sin embargo, en otras áreas, como el disfrute de las matemáticas (Pregunta 3) y el impulso por explorar problemas de manera autónoma (Pregunta 4), los puntajes fueron moderados (3.27 y 3.17, respectivamente), sugiriendo que, aunque hubo interés, la motivación y disfrute no fueron igualmente elevados en todos los aspectos (Figura IX).

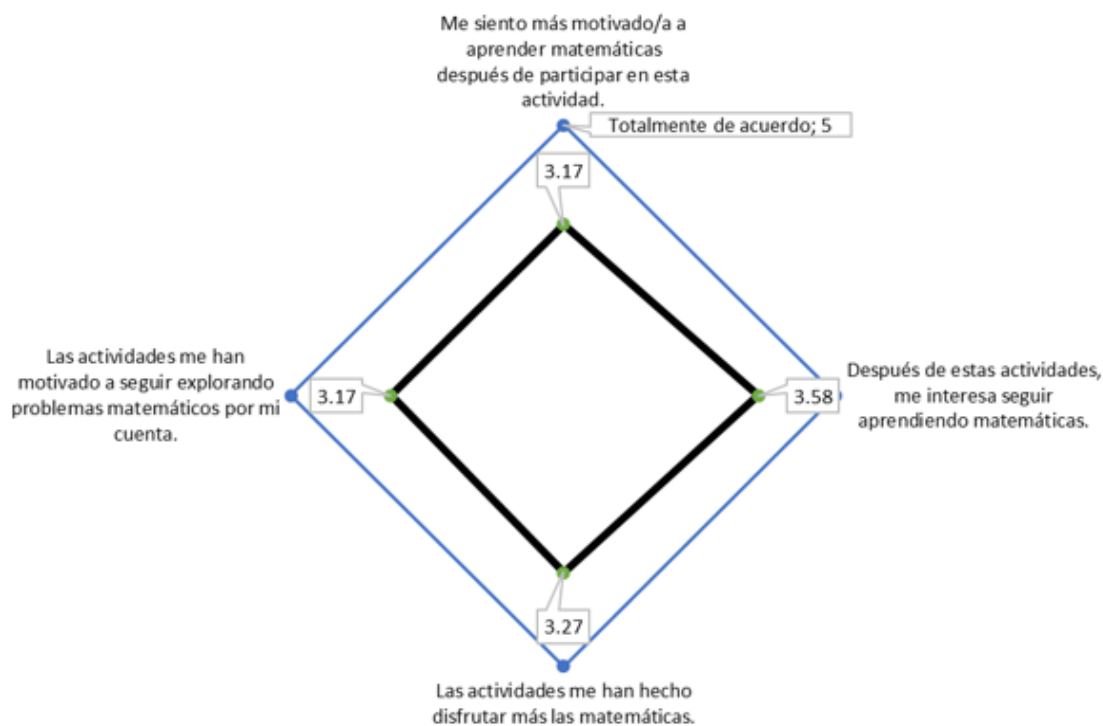

**Figura IX**. Puntaje promedio por pregunta relacionada con la Motivación

En cuanto a la distribución por niveles de motivación, la mayoría de los estudiantes presentó alta motivación en la Pregunta 2, con un 50% reportando un nivel alto. Las Preguntas 3 y 4 también mostraron proporciones significativas en el nivel alto (33% y 38%, respectivamente). Sin embargo, las preguntas sobre motivación general (Pregunta 1) y exploración autónoma (Pregunta 4) mostraron porcentajes más altos en los niveles de motivación media o baja (23% y 25%, respectivamente). Esto resalta la necesidad de fortalecer la motivación autónoma y el disfrute de las matemáticas en futuros enfoques (Figura X).

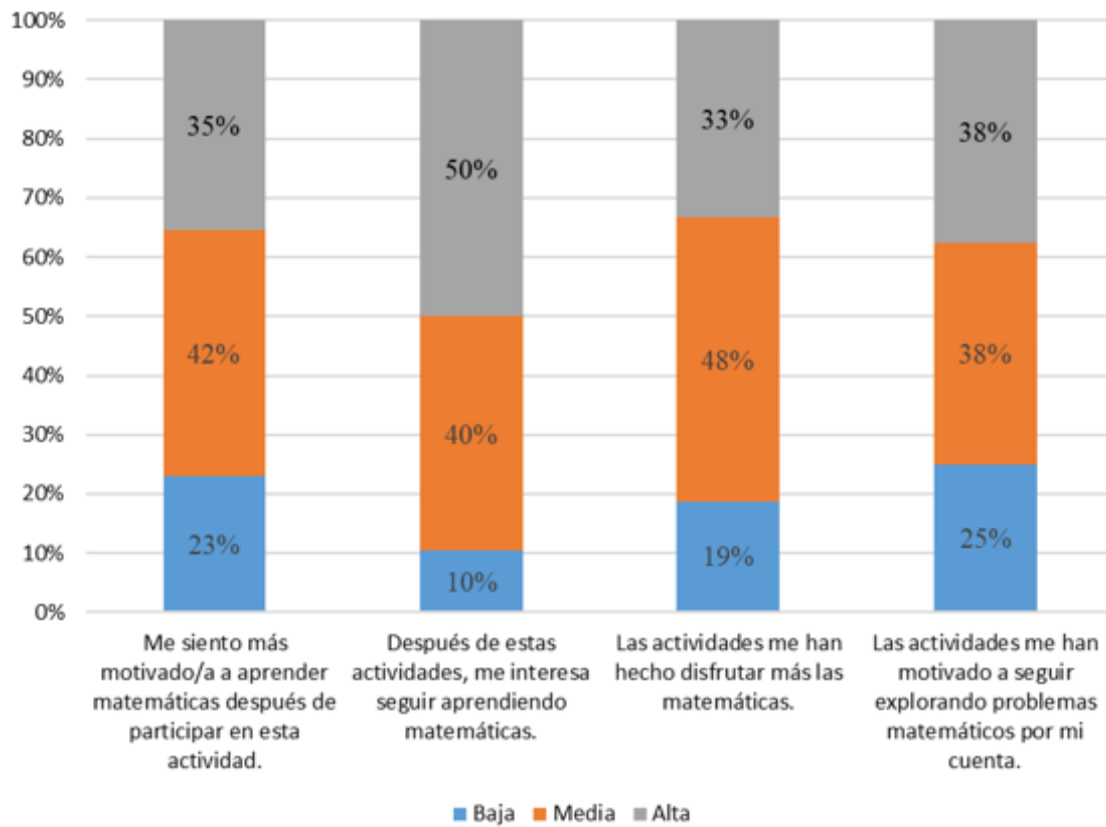

**Figura X**. Niveles de Motivación por pregunta realizada

La relación de estos resultados con los objetivos del estudio es clara: el enfoque Matemáticas en Tres Actos logra generar un interés positivo hacia el aprendizaje matemático, pero aún presenta desafíos en cuanto a la motivación para explorar problemas matemáticos de manera independiente. Para maximizar su impacto, sería importante diseñar estrategias específicas que refuercen el impulso hacia la exploración autónoma y profundicen en el disfrute de las matemáticas, aspectos que mostraron áreas de mejora. Además, la necesidad de intervenir en los estudiantes con niveles de motivación baja en algunas áreas resalta la importancia de crear estrategias que conecten las matemáticas con situaciones cotidianas, haciendo la materia más accesible y atractiva para todos los estudiantes.

**Resultados relacionados con la Actitud hacia las Matemáticas**

Los estudiantes mostraron una mejora general en su actitud hacia las matemáticas, destacando especialmente la percepción de la utilidad de esta disciplina en la vida cotidiana (Pregunta 2), que obtuvo el promedio más alto (3.88). Esto sugiere que el enfoque implementado contribuyó a una mayor apreciación de las matemáticas como una herramienta práctica y relevante. Sin embargo, la percepción de cercanía y comprensión de las matemáticas (Pregunta 4) obtuvo el promedio más bajo (3.08), indicando que algunos estudiantes aún enfrentan desafíos para conectar las matemáticas con experiencias personales.

Los promedios intermedios en las Preguntas 1 (3.27) y 3 (3.69) sugieren que, aunque los estudiantes mostraron un interés creciente, todavía existen oportunidades para fortalecer su percepción de accesibilidad y utilidad (Figura XI).

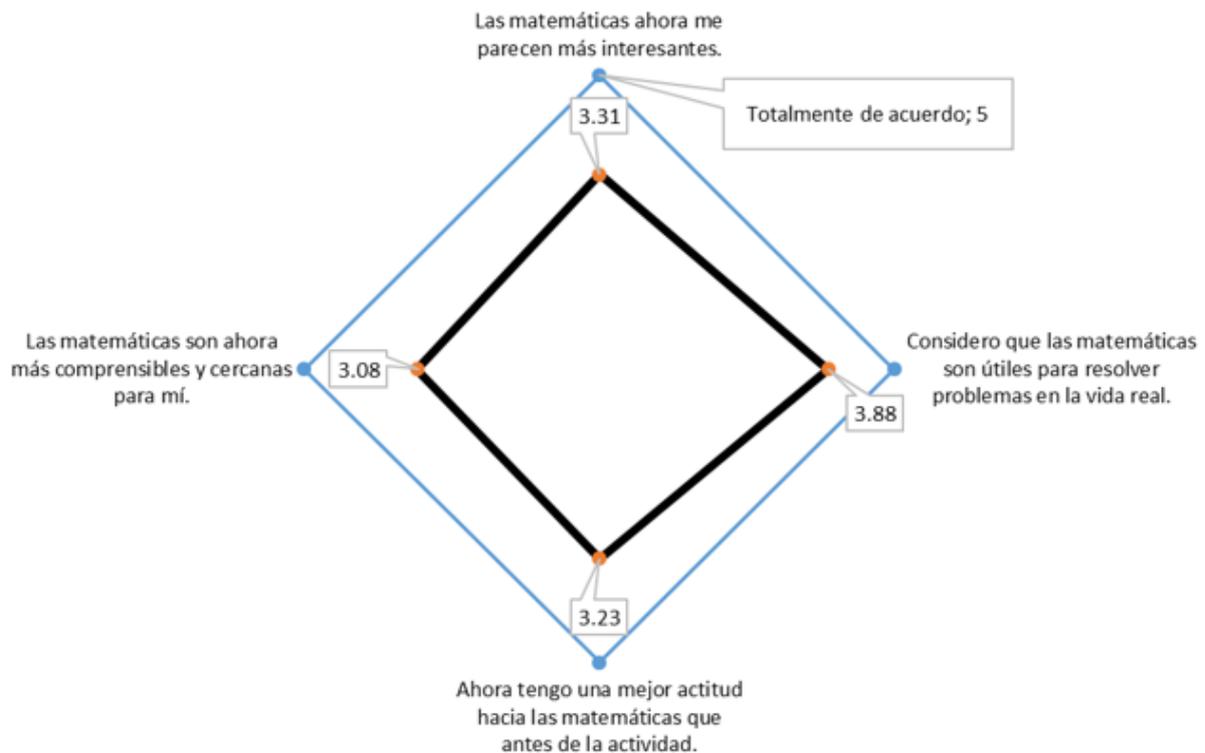

**Figura XI**. Puntaje promedio por pregunta relacionada con la Actitud hacia las Matemáticas

La distribución de respuestas en los niveles de actitud (alta, media y baja) revela patrones consistentes: la alta actitud se destacó principalmente en la Pregunta 2, con un 69% de respuestas, indicando una fuerte valoración de la utilidad práctica de las matemáticas, mientras que solo el 29% de los estudiantes calificaron en este nivel para la Pregunta 4, lo que sugiere una menor percepción de accesibilidad y comprensión. En cuanto a la actitud media, predominó en las Preguntas 3 (44%) y 4 (40%), reflejando que, aunque muchos estudiantes experimentaron mejoras, varios aún mantienen una actitud moderada hacia las matemáticas. Finalmente, la baja actitud se observó en la Pregunta 4, con un 27% de respuestas, lo que destaca la necesidad de implementar estrategias adicionales para hacer las matemáticas más comprensibles y accesibles (Figura XII).

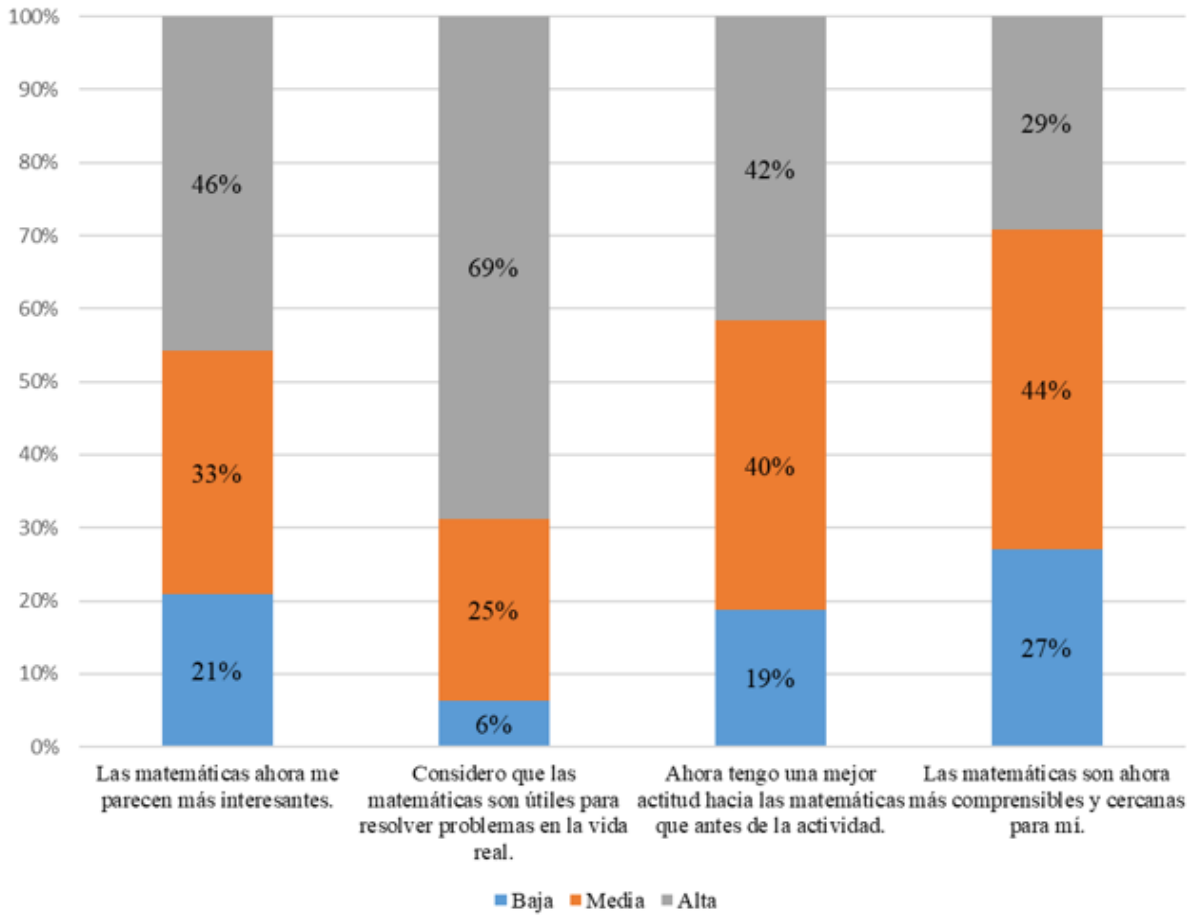

**Figura XII**. Niveles de Actitud hacia las Matemáticas por pregunta realizada

Los resultados indican que el enfoque Matemáticas en Tres Actos tuvo un impacto positivo en la actitud de los estudiantes, particularmente en la percepción de la utilidad práctica de las matemáticas. Este hallazgo es consistente con los objetivos del estudio, que buscan generar una experiencia educativa significativa y atractiva. No obstante, los datos también apuntan a la necesidad de ajustar las estrategias para mejorar la percepción de cercanía y comprensión, especialmente en contextos más prácticos y cotidianos

**Resultados relacionados con la Confianza en su Capacidad Matemática**

En términos generales, los estudiantes reportaron un aumento en su confianza matemática, con un impacto particularmente positivo en su capacidad para interpretar y resolver problemas matemáticos (Pregunta 3), que obtuvo el promedio más alto (3.69). Este hallazgo sugiere que el enfoque promovió una autoconfianza sólida en habilidades prácticas y específicas. Sin embargo, el promedio más bajo (3.10) corresponde a la Pregunta 2, relacionada con la percepción de seguridad al resolver problemas matemáticos. Esto indica que algunos estudiantes experimentaron inseguridad en este aspecto, lo que evidencia la necesidad de diseñar estrategias específicas que refuercen su confianza en situaciones de resolución de problemas complejos.

Los promedios de las Preguntas 1 y 4 fueron similares (3.27 y 3.30, respectivamente), lo que refleja una confianza moderada en la comprensión de conceptos matemáticos y en la habilidad general para abordar problemas matemáticos (Figura XIII)

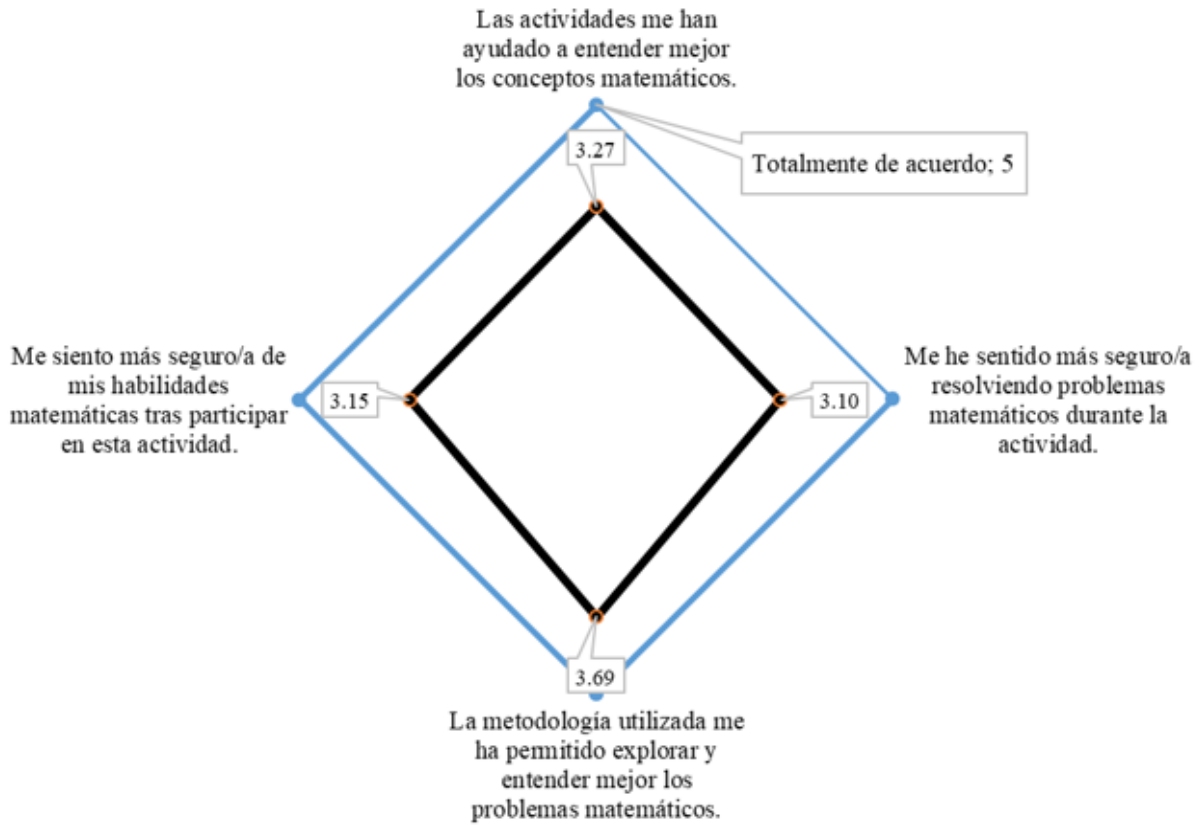

**Figura XIII**. Puntaje promedio por pregunta relacionada con la Confianza en su Capacidad Matemática

La distribución de las respuestas en los niveles de confianza (alta, media y baja) revela patrones interesantes: en alta confianza, la Pregunta 3 lidera con un 52%, indicando que la mayoría de los estudiantes se sintió fortalecida en habilidades prácticas, mientras que las Preguntas 1 y 4 también mostraron proporciones considerables en este nivel (35% y 40%, respectivamente), reflejando un impacto positivo en la confianza matemática general. En cuanto a confianza media, predominó en las Preguntas 1 y 4, con un 48% y 35% de respuestas, lo que sugiere que muchos estudiantes percibieron mejoras moderadas, pero aún no alcanzaron una seguridad plena en sus habilidades. Finalmente, en baja confianza, la mayor proporción de respuestas se observó en la Pregunta 2, con un 27%, lo que subraya que algunos estudiantes continúan enfrentando dificultades para percibir mejoras en su confianza al resolver problemas matemáticos (Figura XIV).

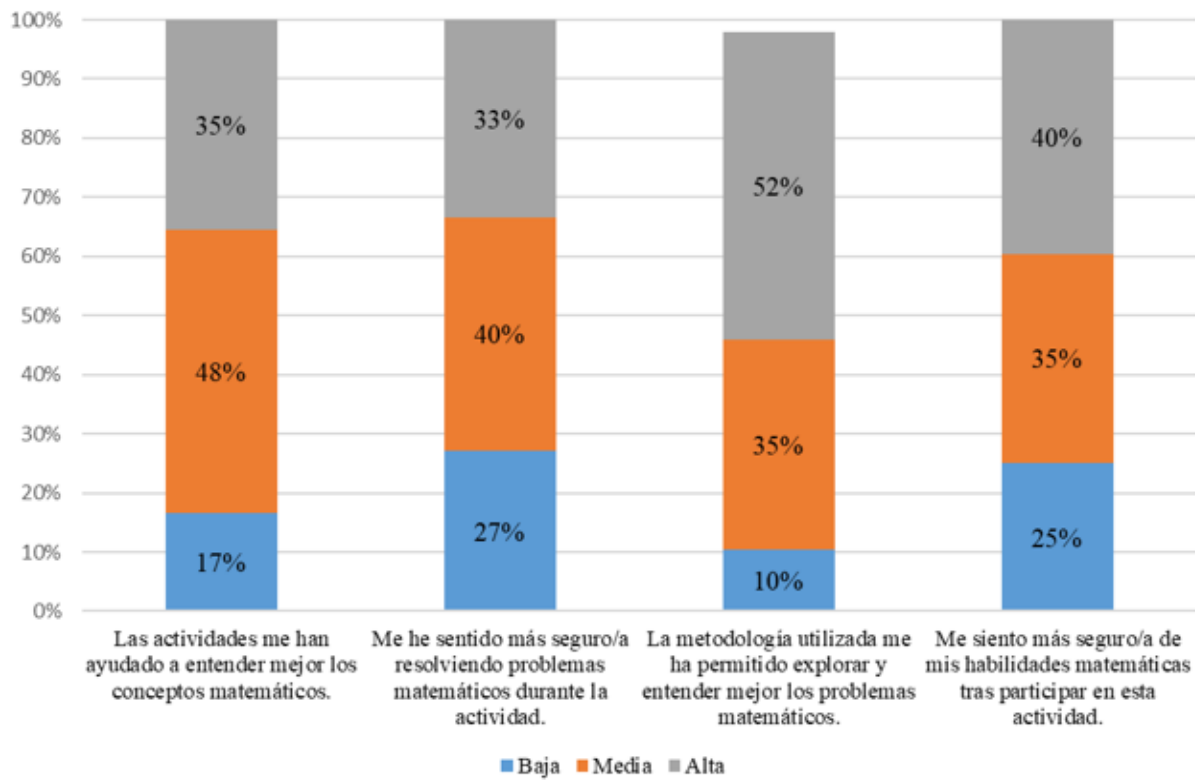

**Figura XIV**. Niveles de Confianza en su Capacidad Matemática por pregunta realizada

Los resultados reflejan que el enfoque Matemáticas en Tres Actos tiene un impacto positivo general en la confianza de los estudiantes hacia sus habilidades matemáticas. Este hallazgo está alineado con los objetivos del estudio, que buscan fomentar una experiencia educativa significativa, donde los estudiantes se sientan empoderados para abordar problemas matemáticos con mayor seguridad. Sin embargo, los datos también subrayan la necesidad de ajustar la metodología para trabajar con mayor énfasis en la resolución de problemas matemáticos, un área donde algunos estudiantes mostraron inseguridad. Actividades que promuevan un refuerzo positivo en habilidades complejas podrían ser clave para mejorar este aspecto

# 5. Conclusiones, Discusión y Recomendaciones

Interpretar y analizamos los hallazgos del estudio en relación con los objetivos específicos establecidos, contrastándolos con la literatura existente y reflexionando sobre sus implicaciones educativas. Además, se abordan las limitaciones del estudio y se presentan recomendaciones para futuras investigaciones en el campo de la educación matemática y la implementación del enfoque Matemáticas en Tres Actos.

# 5.1 Logro del Objetivo General y los Objetivos Específicos

Este estudio ha demostrado que la adaptación del enfoque Matemáticas en Tres Actos tiene un impacto positivo en diversos aspectos del aprendizaje matemático de los estudiantes de secundaria. A través de un análisis detallado, se evaluaron las competencias en creatividad matemática, resolución de problemas, habilidades metacognitivas y percepciones hacia las matemáticas, y se han alcanzado los siguientes resultados:

- Creatividad Matemática: Los estudiantes mostraron una notable capacidad para generar y formular problemas matemáticos de manera creativa, especialmente en los grados superiores. El enfoque facilitó la reinterpretación de contextos y la formulación de preguntas originales, destacándose en los indicadores de Diversidad de Enfoques, Fluidez en la Formulación de Preguntas, Interpretación Creativa del Contexto e Innovación en el Enfoque del Problema. En los grados inferiores, se observaron algunas limitaciones en la originalidad de las preguntas, sugiriendo que se debe seguir trabajando en la estimulación de la creatividad desde etapas tempranas.
- Resolución de Problemas: Los estudiantes demostraron competencias sólidas en la comprensión de los problemas y en la traducción de los mismos al lenguaje matemático. En general, lograron aplicar estrategias matemáticas adecuadas. Sin embargo, en la justificación de procedimientos, especialmente en los grados inferiores, se identificaron áreas de mejora. Los resultados confirman que el enfoque mejora la resolución de problemas y fomenta el razonamiento matemático, aunque es necesario reforzar la argumentación.
- Habilidades Metacognitivas: Los estudiantes mostraron avances significativos en habilidades metacognitivas como la autonomía en la toma de decisiones y la autorregulación. Aunque la flexibilidad cognitiva y la corrección de errores fueron áreas con desafíos, se observó una mejora en la capacidad de autorregulación y reflexión sobre sus propios procesos de pensamiento, lo que es fundamental para un aprendizaje autónomo y efectivo.
- Motivación y Actitudes hacia las Matemáticas: Los estudiantes experimentaron un aumento en su motivación y actitud hacia las matemáticas, mostrando un mayor interés por la asignatura y una mayor confianza en sus habilidades matemáticas. El enfoque contribuyó significativamente a generar una experiencia de aprendizaje más atractiva, lo que a su vez incrementó el compromiso de los estudiantes con las matemáticas.

# 2.2 Implicaciones para la Educación Matemática

Los resultados del estudio subrayan la eficacia del enfoque Matemáticas en Tres Actos para promover una experiencia educativa activa y centrada en el estudiante. Este enfoque no solo mejora las competencias matemáticas, sino que también favorece el desarrollo de habilidades metacognitivas, creatividad y una actitud positiva hacia las matemáticas. Además, es relevante para los docentes que deseen integrar metodologías innovadoras en su aula. Se recomienda:

- Para docentes: Fomentar la creatividad y la resolución de problemas, proporcionando un entorno de aprendizaje que valore las ideas originales y que permita la exploración de diferentes enfoques. Asimismo, es fundamental fomentar el razonamiento y la argumentación matemática a través de actividades que inviten a la reflexión y la justificación de soluciones.
- Para instituciones educativas: Considerar la integración de enfoques activos como el Matemáticas en Tres Actos en sus programas de matemáticas, apoyando a los docentes con recursos y formación adecuada.
- Para políticas educativas: Impulsar la inclusión de metodologías activas en los lineamientos curriculares y promover la formación continua de los docentes en enfoques innovadores.

## 5.3 Limitaciones del Estudio y Recomendaciones para Futuras Investigaciones

Aunque los hallazgos del estudio son prometedores, es importante considerar algunas limitaciones. El tamaño y la representatividad de la muestra son limitados, lo que puede restringir la generalización de los resultados. Además, el estudio se basó en métodos cualitativos, lo que introdujo cierto grado de subjetividad en la interpretación de los resultados. Para futuras investigaciones, se recomienda:

- Ampliar la muestra: Incluir más estudiantes y contextos educativos para obtener resultados más generalizables.
- Estudios longitudinales: Realizar investigaciones a largo plazo para evaluar el impacto sostenido del enfoque en las competencias matemáticas y las habilidades metacognitivas.
- Enfoques mixtos: Incorporar métodos cuantitativos que permitan obtener una visión más objetiva y robusta sobre los efectos del método.

## 5.4 Reflexión Final

En conclusión, el enfoque Matemáticas en Tres Actos ha demostrado ser una metodología eficaz para mejorar las habilidades matemáticas, metacognitivas y afectivas de los estudiantes de secundaria. La creatividad matemática, el desarrollo de competencias en resolución de problemas, las habilidades metacognitivas y la motivación hacia las matemáticas se han visto significativamente mejoradas, lo que confirma la validez de este enfoque para la educación matemática. La implementación de métodos centrados en el estudiante que conecten las matemáticas con experiencias cotidianas y fomenten la autorregulación es fundamental para preparar a los estudiantes para enfrentar los retos académicos y personales de manera integral.